\RequirePackage{lineno}
\documentclass[11pt]{article}
\usepackage[utf8]{inputenc}
\usepackage[english]{babel}
\usepackage{graphicx}
\usepackage{amsmath}
\usepackage{amssymb}
\usepackage{amsfonts}
\usepackage{amsthm}
\usepackage[left=2.1cm,top=1.5cm,right=2.1cm, bottom=2.2cm,letterpaper]{geometry}
\usepackage{mathrsfs}
\usepackage{caption}
\usepackage{subcaption}
\usepackage{latexsym}
\usepackage{xfrac}
\usepackage{tcolorbox}
\usepackage{enumitem}
\usepackage[symbol]{footmisc}
\usepackage{float}
\usepackage{hyperref}
 \usepackage{url}
\usepackage[toc]{appendix}
\newtheorem{theorem}{Theorem}[section]

\newtheorem{proposition}{Proposition}[section]

\newtheorem{remark}{Remark}[section]
\numberwithin{equation}{section}
\numberwithin{figure}{section}

\newcommand{\R}{\mathbb{R}}

\makeatletter
\renewcommand*\env@matrix[1][*\c@MaxMatrixCols c]{%
	\hskip -\arraycolsep
	\let\@ifnextchar\new@ifnextchar
	\array{#1}}
\makeatother

\def\proof{{\bf Proof.}\ }

\def\eq#1{(\ref{#1})}
\def\neweq#1{\begin{equation}\label{#1}}
\def\endeq{\end{equation}}

\begin{document}

\title{Solenoidal extensions in domains with obstacles:\\
explicit bounds and applications to Navier-Stokes equations}

\author{Ilaria FRAGAL\`{A} - Filippo GAZZOLA - Gianmarco SPERONE\footnote{Corresponding author. E-mail address: \textbf{sperone@karlin.mff.cuni.cz}}}
\date{}
\maketitle
\vspace*{-6mm}

\begin{abstract} We introduce a new method for constructing solenoidal extensions of fairly general boundary data
in (2d or 3d) cubes that contain an obstacle. This method allows us  to provide {\it explicit} bounds for the Dirichlet norm of the extensions.
It runs as follows: by inverting the trace operator, we first determine suitable extensions, not necessarily solenoidal, of the data; then we analyze the Bogovskii problem with the resulting divergence to obtain a solenoidal extension; finally, by solving a variational problem involving the infinity-Laplacian and using ad hoc cutoff functions, we find explicit bounds in terms of the geometric parameters of the obstacle. The natural applications of our results lie in the analysis of inflow-outflow problems, in which an explicit bound on the inflow velocity is needed to estimate the threshold for uniqueness in the stationary Navier-Stokes equations and, in case of symmetry, the stability of the obstacle
immersed in the fluid flow.\par
\end{abstract}

\section{Introduction}

\begin{figure}[H]
	\begin{center}
		\scalebox{-1}[1]{\includegraphics[height=55mm]{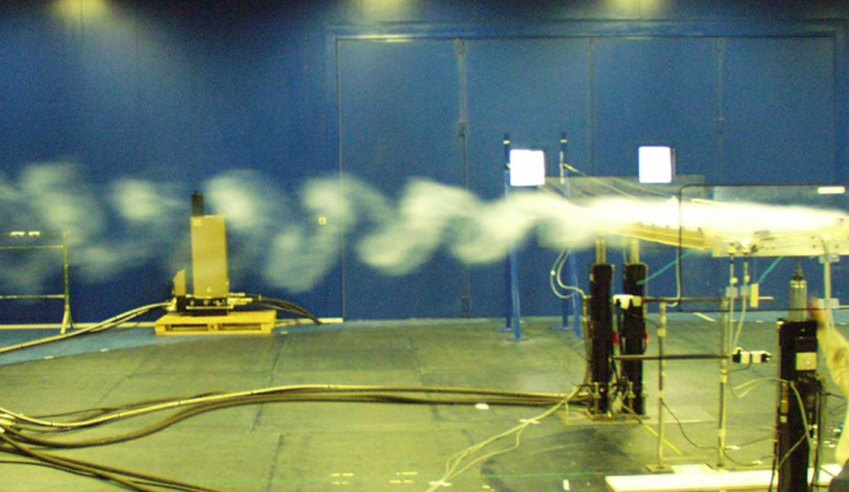}}		
	\end{center}
	\vspace*{-5mm}
	\caption{Vortices around a plate obtained in wind tunnel experiments at Politecnico di Milano.}\label{WT}
\end{figure}

Stationary inflow-outflow problems in fluid mechanics are well-modeled by the (steady state) Navier-Stokes equations
describing the motion of the fluid and by {\em nonhomogeneous} boundary conditions prescribing how a given fluid enters or exits the
considered bounded domain $\Omega$ (either in $\R^2$ or in $\R^3$):
\begin{equation}\label{nsstokes0}
-\eta\Delta u+(u\cdot\nabla)u+\nabla p=0,\ \quad  \nabla\cdot u=0 \ \mbox{ in } \ \Omega,\qquad u=h \ \mbox{ on } \ \partial\Omega\,.
\end{equation}
In \eqref{nsstokes0}, $u$ is the velocity vector field, $p$ is the
scalar pressure and $\eta>0$ is the kinematic viscosity, while
the datum $h$ contains both the inflow-outflow conditions and the behavior on the remaining part of the boundary. The theory  developed so far
in order to manage Navier-Stokes equations under nonhomogeneous
Dirichlet conditions
(see e.g.\ \cite{galdi2011introduction})
suggests to reduce the problem to homogeneous conditions through a suitable {\em solenoidal extension} of the boundary data. Namely, one needs to find a vector field $v_0$ satisfying
\begin{equation}\label{div}
\nabla\cdot v_0=0 \ \text{ in } \ \Omega,\qquad v_0=h \ \text{ on } \ \partial\Omega\, .
\end{equation}
This problem, whose interest and applicability go far beyond fluid mechanics, has a long history, starting from the pioneering works
of Sobolev \cite{sobolev1954new}, Cattabriga \cite{cattabriga1961problema} and Ladyzhenskaya-Solonnikov \cite{ladyzhenskaya1969mathematical,ladyzhenskaya1978some}; see also the book by Galdi \cite[Section III.3]{galdi2011introduction}.
If the boundary conditions are themselves of solenoidal type such as constants, Poiseuille or Couette flows
(see \cite{landau} for other models), the extension is found by a fairly standard procedure, see \cite{lith1,solonn1,solonn2} for bounded domains and \cite{chipot1,chipot2} for a special class of unbounded domains.
The classical way to solve \eqref{div} relies in the use of a proper extension of the data $h$ as a curl,
together with a Hopf's-type cutoff function, see \cite[p.130]{ladyzhenskaya1969mathematical} and also \cite[Section IX.4]{galdi2011introduction}.
However, if the inflow-outflow datum $h$ does not have a straightforward solenoidal extension, the problem
becomes significantly more difficult. This is the case, for instance, when the considered domain contains an {\em obstacle} where, due to the effects of viscosity, the flow satisfies no-slip conditions.
Then, even if the inflow-outflow datum has a simple solenoidal
extension, one can still use cutoff functions in order to meet the homogeneous boundary conditions on the obstacle. Nevertheless, it turns out that,
in real life, the obstacle perturbs the fluid flow, creating some {\em vortices} behind itself even for low-Reynolds-numbers. Recent experimental
and numerical evidence \cite{bordogna2020effects,uruba2019reynolds,uruba3d2} shows that, at a sufficiently large distance behind the obstacle,
the perturbed flow concentrates its turbulent motion mostly in the wake of the body, see Figure \ref{WT} for a wind tunnel experiment and \cite{bordogna2020effects,rocchi} for some experimental data. 
 
It is then clear that different boundary conditions should be imposed on the outlet, see \cite{gazzolasperone1,heywood1996artificial,kravcmar2018modeling}. Overall, in presence of an obstacle, the most realistic physical boundary conditions
are different (nonhomogeneous) inflow and outflow conditions, combined with homogenous conditions on the obstacle. In this situation,
the construction of a solenoidal extension appears possible only in two steps. First, to find an extension, not necessarily solenoidal, of the
inflow-outflow data, thereby ``inverting'' the trace operator for vector fields; this problem has been systematically studied since the works of Miranda \cite{miranda1951sulla}, Prodi \cite{prodi1956tracce} and Gagliardo \cite{gagliardo}, giving a full characterization of the trace operator and providing an explicit extension for any locally Lipschitz domain and any boundary datum. Second, to solve the Bogovskii problem \cite{bogovskii1979solution,bogovskii1980} with the resulting divergence: the celebrated Bogovskii formula, dating back to 1979, yields a class
of solutions by means of the Calderón-Zygmund theory of singular integrals. Durán \cite{duran2012elementary} proposed in 2012 an alternative approach based on the Fourier transform. Incidentally, let us also mention that the Bogovskii problem is strictly related to several
inequalities arising both in fluid mechanics and elasticity; see \cite{acosta2017divergence,babuska,costabel,friedrichs,horgan,korn}.

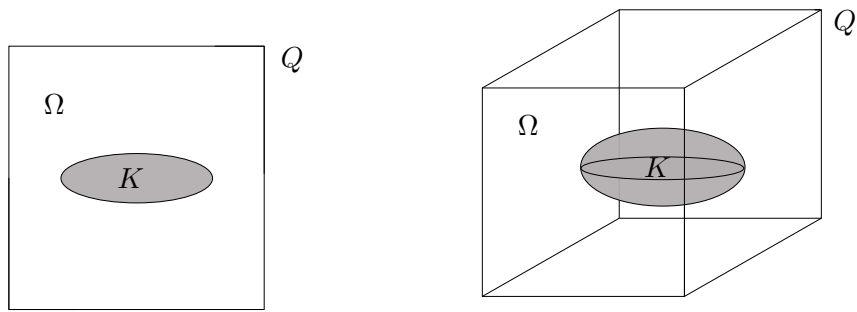
\begin{figure}[H]
	\null\vskip - 5.5 cm
	\begin{minipage}{0.5\linewidth}
		\centering
		\def\svgwidth{7cm}
		\begingroup%
		\makeatletter%
		\providecommand\color[2][]{%
			\errmessage{(Inkscape) Color is used for the text in Inkscape, but the package 'color.sty' is not loaded}%
			\renewcommand\color[2][]{}%
		}%
		\providecommand\transparent[1]{%
			\errmessage{(Inkscape) Transparency is used (non-zero) for the text in Inkscape, but the package 'transparent.sty' is not loaded}%
			\renewcommand\transparent[1]{}%
		}%
		\providecommand\rotatebox[2]{#2}%
		\ifx\svgwidth\undefined%
		\setlength{\unitlength}{377.58676bp}%
		\ifx\svgscale\undefined%
		\relax%
		\else%
		\setlength{\unitlength}{\unitlength * \real{\svgscale}}%
		\fi%
		\else%
		\setlength{\unitlength}{\svgwidth}%
		\fi%
		\global\let\svgwidth\undefined%
		\global\let\svgscale\undefined%
		\makeatother%
		\begin{picture}(1,1.38846134)%
		\put(0,0)
		{\hskip 2.3 cm \includegraphics[height=40mm]{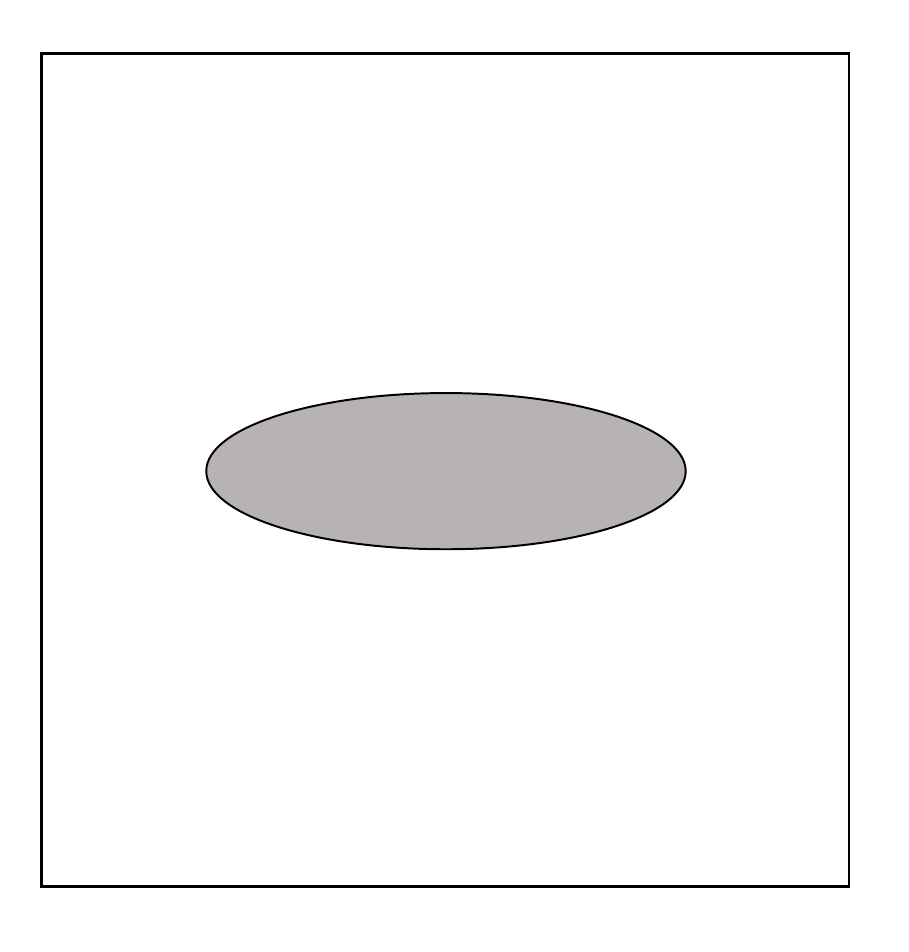}
			\hskip 2.3 cm
			\includegraphics[height=45mm]{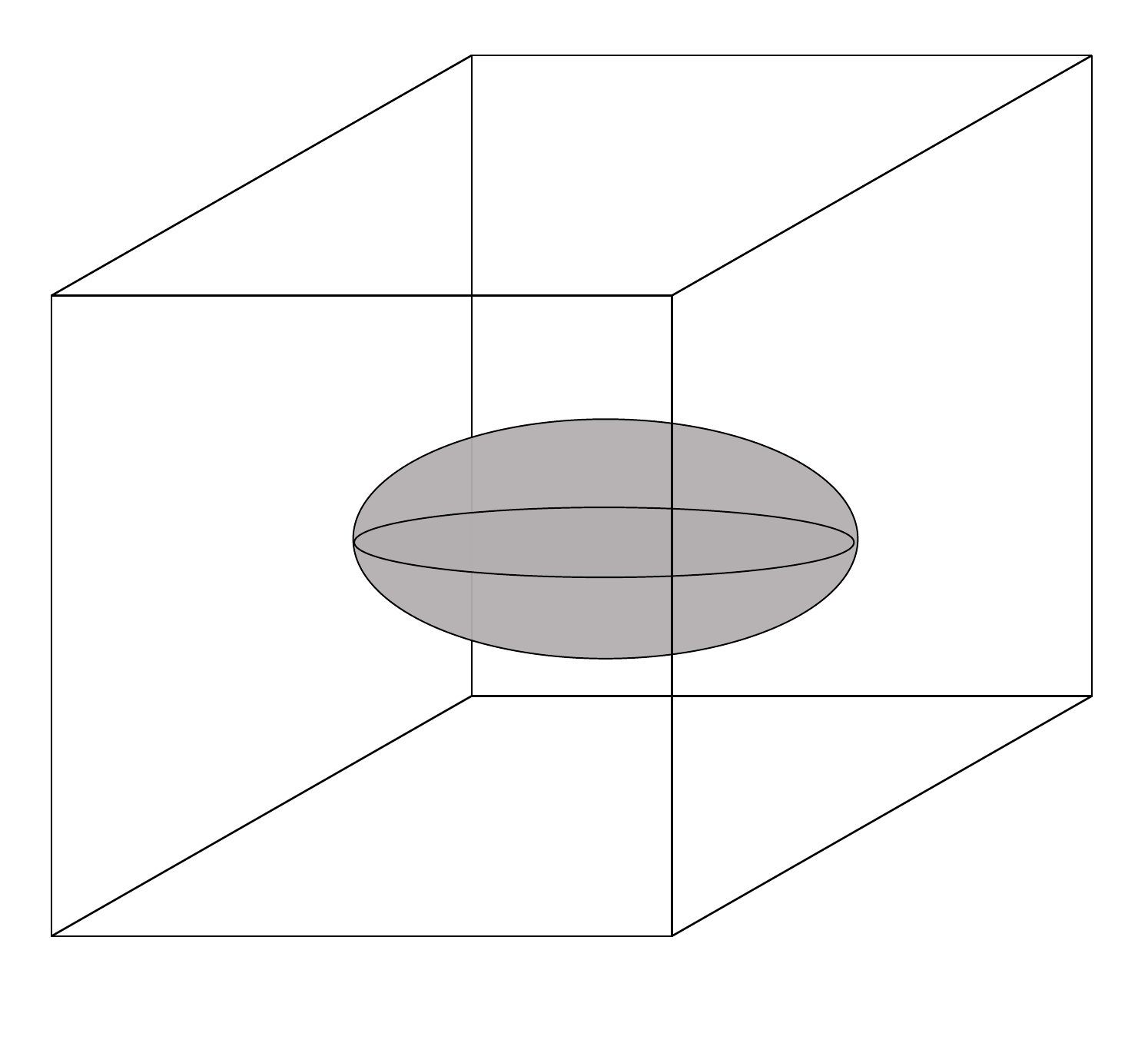}}
		\put(0.87,0.5){\color[rgb]{0,0,0}\makebox(0,0)[lb]{\smash{$Q$}}}%
		\put(0.56,0.27){\color[rgb]{0,0,0}\makebox(0,0)[lb]{\smash{$K$}}}%
		\put(0.42,0.41){\color[rgb]{0,0,0}\makebox(0,0)[lb]{\smash{$\Omega$}}}%
		\put(1.92,0.57){\color[rgb]{0,0,0}\makebox(0,0)[lb]{\smash{$Q$}}}%
		\put(1.56,0.29){\color[rgb]{0,0,0}\makebox(0,0)[lb]{\smash{$K$}}}%
		\put(1.32,0.37){\color[rgb]{0,0,0}\makebox(0,0)[lb]{\smash{$\Omega$}}}%
		\end{picture}%
		\endgroup%
	\end{minipage}%
	\begin{minipage}{0.5\linewidth}
		\centering
		\def\svgwidth{7cm}
		\begingroup%
		\makeatletter%
		\providecommand\color[2][]{%
			\errmessage{(Inkscape) Color is used for the text in Inkscape, but the package 'color.sty' is not loaded}%
			\renewcommand\color[2][]{}%
		}%
		\providecommand\transparent[1]{%
			\errmessage{(Inkscape) Transparency is used (non-zero) for the text in Inkscape, but the package 'transparent.sty' is not loaded}%
			\renewcommand\transparent[1]{}%
		}%
		\providecommand\rotatebox[2]{#2}%
		\ifx\svgwidth\undefined%
		\setlength{\unitlength}{379.39037574bp}%
		\ifx\svgscale\undefined%
		\relax%
		\else%
		\setlength{\unitlength}{\unitlength * \real{\svgscale}}%
		\fi%
		\else%
		\setlength{\unitlength}{\svgwidth}%
		\fi%
		\global\let\svgwidth\undefined%
		\global\let\svgscale\undefined%
		\makeatother%
		\endgroup%
	\end{minipage}
	\caption{The obstacle $K$ in the box $Q$,  in $2d$  (left) and in $3d$ (right).}
	\label{sheet}
\end{figure}

As a consequence of the vortex shedding, the fluid exerts forces on the obstacle and, if one is interested in the stability of the obstacle itself, the most relevant one is the lift force. With the final purpose of analyzing the stability of a suspension bridge under the action of the wind \cite{gazzolabook}, a simplified geometric framework where to analyze the appearance a of lift force was first suggested in \cite{phd} and subsequently discussed in \cite{bongazspe,gazzolasperone1,gazspe}.
The setting can be simply described as follows: the container $Q$ is given by an open square in $\R ^2$ or an open cube in $\R^3$,
while the obstacle is represented as a compact, connected, and simply connected domain $K$ with Lipschitz boundary contained into $Q$, see Figure \ref{sheet}. In this geometric setting, the main purpose of this paper is the construction of estimable solenoidal extensions of quite
general boundary data to $\Omega=Q \setminus K$. Precisely, in dimension $d=2$ or $d=3$, given a vector field $h$ satisfying suitable assumptions, we determine a solution $v_0$ to the boundary value problem \eqref{div}, along with some upper bound on the Dirichlet norm of $v_0$ in $\Omega$. The goal is not simply to show the existence of some solenoidal extension, but also to obtain an {\em explicit form} of it, in order to derive  {\it explicit bounds} on its norm. The reason is that we are mainly interested in applications to fluid mechanics, such as finding bounds on the inflow velocity guaranteeing the unique solvability of the Navier-Stokes equations \eqref{nsstokes0}. In turn, in view of the results contained in \cite{gazspe}, in a symmetric framework, unique solvability implies that the lift applied over $K$ is zero.\par
The paper is organized as follows. Section \ref{glance} serves as a guideline, where the outline of our strategy is presented, together with our main result (Theorem \ref{teononconstant}) and its application to Navier-Stokes equations. The steps of this strategy are carried over in the remaining sections of the article. In Section \ref{functional} we formulate an extension result, see Theorem \ref{gadliardo}, that not only allows to invert the trace operator in our geometric framework, but also to study two new inflow-outflow models that are suggested. Then, in Section \ref{sec:phi} we explicitly solve a variational problem involving the infinity-Laplacian, yielding a sharp bound on the $W^{1,\infty}$-norm of a class of scalar cutoff functions. Through a delicate combination of the results contained in \cite{bogovskii1979solution, bogovskii1980, duran2012elementary,galdi2011introduction}, an upper bound for the Bogovskii constant of the domain $\Omega$ is found in Section \ref{duransec}, see Theorem \ref{unionstarshapedcor}; this requires the estimation of the norms of certain mollifiers given in Section \ref{normsmollifiers}. By using these results we finally give the proof of Theorem \ref{teononconstant} in Section \ref{finalproof}.

\section{The paper at a glance}\label{glance}

\subsection{Assumptions and outline of the strategy}\label{outline}

We let $Q=(-L,L)^{d}$ and $K \subset Q \subset \R ^ d$ ($d = 2$ or $d= 3$) be as described above, see again Figure \ref{sheet}. We fix here the assumptions  on the boundary datum $h$ in problem \eqref{div}:
we view $h$ as the compound of four fields $h_{i}$ defined on the four sides  $\Sigma_i$'s of $Q$ in dimension $2$, and of
six fields $h_{i}$ defined on the six faces $\Sigma_i$'s of $Q$ in dimension $3$.
We denote by $\hat{n}$ the outward unit normal to (the sides/faces of) $Q$ and by $\mathcal V$ the family of vertices of $Q$:

\begin{equation} \label{zeroflux2}
h \in \mathcal C (\partial Q)\, , \qquad h _i \in H ^ 1 (\Sigma_i) \,
, \qquad \int_{\partial Q} h \cdot \hat{n} = 0\, , \qquad h{ \big | _\mathcal V } = 0 \ \ \text{ if } \ \ d = 3\,.
\end{equation}
Note that the vanishing condition for $h$ at the vertices of $Q$ is assumed only for $d = 3$, while it is not needed
for the validity of our results in dimension $d = 2$.
We fix real numbers $a, b, c$ so that
\begin{equation}\label{OmegaR}
L >  a \geq b \geq c >0  \quad \text{ and } \quad K \subset  P  \subset Q \, , \quad \text{ with } \ \begin{cases}
P=(-a,a)\times(-b,b) & \ \text{ if $d=2$}
\\  P=(-a,a)\times(-b,b) \times (-c, c) & \ \text{ if $d=3$} \, ,
\end{cases}
\end{equation}so that the obstacle $K$ is enclosed by the parallelepiped $P$. Throughout the paper we set
\neweq{OmegaRR}
\Omega_{0} \doteq Q\setminus \overline P\,.
\endeq
Working on the set $\Omega_{0}$ allows us to obtain explicit bounds; with our approach, it is clear that the best possible bounds are found by taking $P$ as the smallest parallelepiped enclosing $K$.
We shall proceed in four steps:

\medskip
{\bf Step 1.} We determine a  vector field $$
A_1 \in H^{1}(Q) \cap \mathcal{C}(\overline{Q}) \quad \mbox{with} \quad A_1=h \ \ \mbox{on} \ \ \partial Q\, ,
$$
such that the $H ^ 1$-norm of $A _1$ on $Q$ can be explicitly computed in terms of the $H ^ 1$-norms of the fields $h _ i$ on $\Sigma_i$. The expression of this (not necessarily solenoidal) field is given in Theorem \ref{gadliardo}.

\medskip
{\bf Step 2.} We construct a scalar function $\phi \in W^{1,\infty}(Q)$ such that
$$
\phi = 1 \ \ \text{on} \ \ \partial Q, \qquad \phi = 0 \ \ \text{in} \ \ \overline P\, ,
$$
in such a way that the  $W ^{1, \infty}$-norm of $\phi$ in $\Omega_{0}$ is explicitly computable and as small as possible.
This function $\phi$ is defined in Section \ref{sec:phi}, and its minimizing property  (which is obtained by using the variational properties of infinity-harmonic functions) is stated in
Theorem \ref{thm:infty}; incidentally, we point out that this bound is sharp.
Then the vector field  $A_{2} \doteq \phi A_{1}$ satisfies
$$
A_{2}=h \ \ \mbox{on} \ \ \partial Q \, , \qquad A_2 = 0  \ \ \text{ in } \ \  \overline P\,.
$$
By using the $H ^1$-bounds on $A_2$ and the $W ^ {1, \infty}$-bounds on $\phi$, we can control the $H ^1$-norm of $A_2$ on $Q$.

\begin{figure}[H]
\begin{minipage}{0.5\linewidth}
\centering
\def\svgwidth{7cm}
\begingroup%
  \makeatletter%
  \providecommand\color[2][]{%
    \errmessage{(Inkscape) Color is used for the text in Inkscape, but the package 'color.sty' is not loaded}%
    \renewcommand\color[2][]{}%
  }%
  \providecommand\transparent[1]{%
    \errmessage{(Inkscape) Transparency is used (non-zero) for the text in Inkscape, but the package 'transparent.sty' is not loaded}%
    \renewcommand\transparent[1]{}%
  }%
  \providecommand\rotatebox[2]{#2}%
  \ifx\svgwidth\undefined%
    \setlength{\unitlength}{377.58676bp}%
    \ifx\svgscale\undefined%
      \relax%
    \else%
      \setlength{\unitlength}{\unitlength * \real{\svgscale}}%
    \fi%
  \else%
    \setlength{\unitlength}{\svgwidth}%
  \fi%
  \global\let\svgwidth\undefined%
  \global\let\svgscale\undefined%
  \makeatother%
  \begin{picture}(1,1.38846134)%
    \put(0,0)
   {\includegraphics[height=55mm]{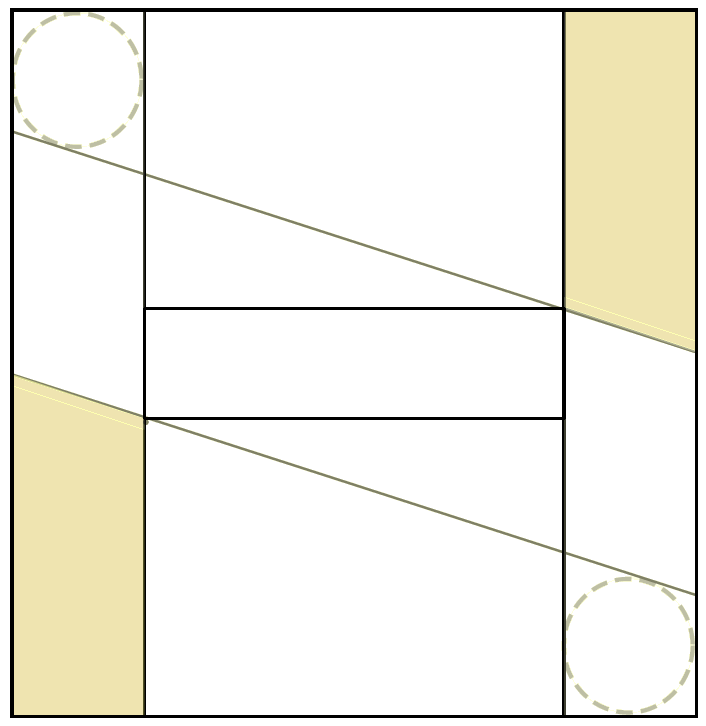}
		\hskip 2 cm
		\includegraphics[height=65mm]{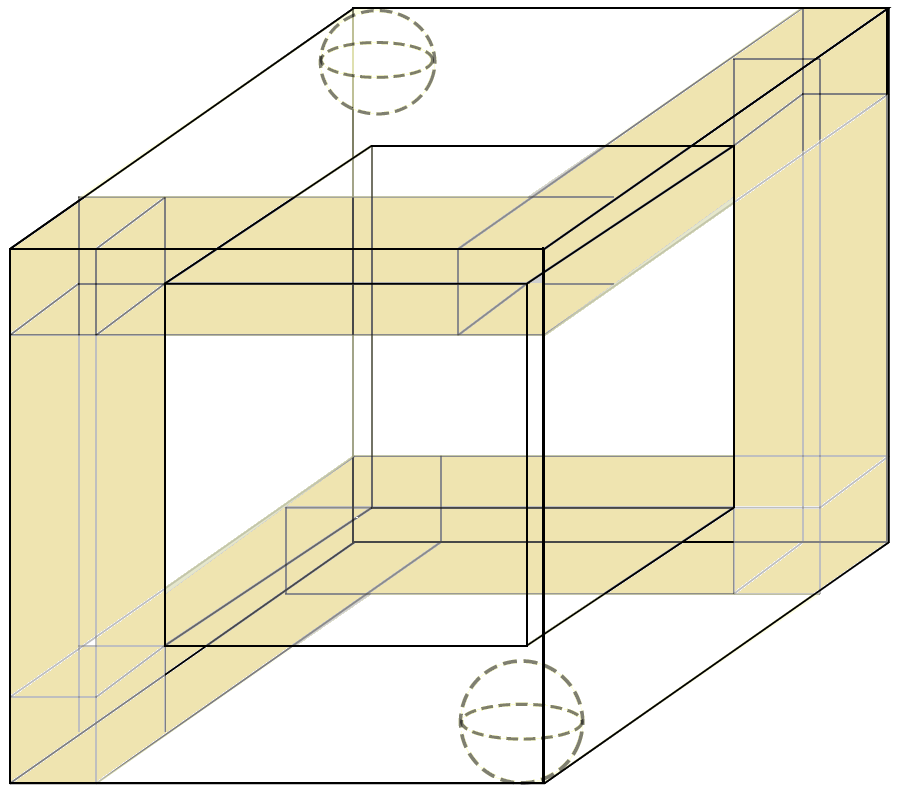}}
  \end{picture}%
\endgroup%
\end{minipage}%
\begin{minipage}{0.5\linewidth}
\centering
\def\svgwidth{7cm}
\begingroup%
  \makeatletter%
  \providecommand\color[2][]{%
    \errmessage{(Inkscape) Color is used for the text in Inkscape, but the package 'color.sty' is not loaded}%
    \renewcommand\color[2][]{}%
  }%
  \providecommand\transparent[1]{%
    \errmessage{(Inkscape) Transparency is used (non-zero) for the text in Inkscape, but the package 'transparent.sty' is not loaded}%
    \renewcommand\transparent[1]{}%
  }%
  \providecommand\rotatebox[2]{#2}%
  \ifx\svgwidth\undefined%
    \setlength{\unitlength}{379.39037574bp}%
    \ifx\svgscale\undefined%
      \relax%
    \else%
      \setlength{\unitlength}{\unitlength * \real{\svgscale}}%
    \fi%
  \else%
    \setlength{\unitlength}{\svgwidth}%
  \fi%
  \global\let\svgwidth\undefined%
  \global\let\svgscale\undefined%
  \makeatother%
  \endgroup%
\end{minipage}
\caption{The (colored) ``doubly-illuminated" region $\Omega_1 \cap \Omega _2$, in 2d (left) and in 3d (right).}
\label{stars}
\end{figure}

\medskip
{\bf Step 3.} We construct a vector field
$$
 A_3\in H^1_0(\Omega_{0}) \ \ \text{satisfying} \ \ \nabla\cdot A_3=-\nabla\cdot A_2 \ \ \mbox{in} \ \ \Omega_{0}\, ,
$$
such that the $L ^ 2$-norm of $\nabla A _ 3$ can be estimated in terms of the $L ^2$-norm of $\nabla  \cdot A_2$.
 To this aim, we provide an upper bound for  the {\it Bogovskii constant of $\Omega_{0}$} (and, hence, of $\Omega$), defined as
\begin{equation}\label{bogo} C_B (\Omega_{0}) \doteq   \sup_ {g \in L ^ 2 _0 (\Omega_{0}) \setminus \{ 0 \} }  \inf  \left\{ \frac{  \|\nabla v\|_{L^2(\Omega_{0})} }{ \|g\|_{L^2(\Omega_{0})}} \ \Bigg| \ v \in H ^ 1 _ 0 (\Omega_{0})  \, , \ \nabla\cdot v=g \ \text{ in } \ \Omega_{0} \right\}  \, ,
\end{equation}
where $L ^ 2 _0 (\Omega_{0})$ denotes the subspace of functions in  $L ^ 2 (\Omega_{0} )$ having zero mean value.
This upper bound is stated Theorem \ref{unionstarshapedcor} and has its own independent interest.  The proof is based on the following  idea:
 first we derive a quantitative version of a result by Durán  \cite{duran2012elementary},
allowing to estimate the Bogovskii constant of a domain which is star-shaped with respect to a ball (see Proposition \ref{p:starshaped}),
and then we apply such result
after decomposing
$\Omega_{0} $ as the union of two domains $\Omega_{1}$ and $\Omega_{2}$ that are star-shaped with respect to a ball placed in a corner: these sets are the regions ``illuminated'' by spherical lamps placed in two opposite corners (namely, each one tangent to the sides of $Q$ intersecting at the corner), see Section \ref{duransec} for the analytic description. In Figure \ref{stars} we illustrate the intersection  $\Omega_1 \cap \Omega _2$ as the colored region ``doubly-illuminated" by both lamps.

%

\vskip -3cm

\medskip
{\bf Step 4 - Conclusion.} Finally,  we observe that the field $v_0$ defined by $v_0 \doteq A_2+A_3$ in $\Omega_{0}$, and extended by zero on $P \setminus K$, is a solution to problem \eqref{div} in $\Omega$; moreover, by making use of  the preceding steps we are in position to find an explicit bound for the $L ^2$-norm of  $\nabla v_0$, see \eqref{f:gb}-\eqref{f:gamma}.

\subsection{Main result and applications to the Navier-Stokes equations}

Putting together all the previously described steps, we now state our main result.

\begin{theorem} \label{teononconstant}
Under assumptions \eqref{zeroflux2}-\eqref{OmegaR}, there exists a vector field $v_0\in H^{1}(\Omega)$ satisfying \eqref{div}
such that
\begin{equation}\label{f:gb}
 \| \nabla v_{0} \|_{L^{2}(\Omega)}  \leq    \Gamma\, ,
 \end{equation}
where $\Gamma$ is a positive constant depending only on $L, a, b, c$ and the norms in $H ^ 1 (\Sigma_i)$ of the functions $h _i$.
Specifically, the value of $\Gamma$ is given by
\begin{equation}\label{f:gamma}\Gamma= \frac{1 + M}{L-a} \,  \| A_{1} \|_{L^{2}( Q)}
+  \| \nabla A_{1} \|_{L^{2}(Q)} + M \, \| \nabla \cdot A_{1} \|_{L^{2}(Q)} \,,
\end{equation}
where $A_1$ is the extension of $h$ found in Step 1, and $M$ is the upper bound for the Bogovskii constant of $\Omega_{0}$ found in Step 3.

%
\end{theorem}

\smallskip

\begin{remark} Concerning the value of the constant $\Gamma$ given in  \eqref{f:gamma}, we further precise that:\par
\noindent
- It is explicitly computable, relying on the explicit expressions of the field $A_1$ and of the constant $M$; to simplify the presentation, such expressions are postponed to Theorems \ref{gadliardo} and \ref{unionstarshapedcor}.\par
\noindent
- It is far from being optimal. In particular, it might be improved by refining the estimates in Steps 1 and 3 (whereas the estimate for the function $\phi$ in Step 2 is sharp).
\end{remark}


\par

Thanks to Theorem \ref{teononconstant} we can give necessary conditions for the appearance of lift forces over an obstacle exerted by Navier-Stokes flows. Let us consider equations \eqref{nsstokes0} with the boundary datum $h$ satisfying \eq{zeroflux2} and a no-slip condition on the obstacle
\begin{equation} \label{boundaryc0}
h=0 \ \text{ on } \ \partial K.
\end{equation}
Since $h\neq0$ on $\partial Q$, we need to deal with both the Sobolev space $H_{0}^{1}(\Omega)$ and the space of functions vanishing only on $\partial K$, which is a proper
connected part of $\partial\Omega$:
$$
H^1_*(\Omega)=\{v\in H^1(\Omega) \ | \ v=0 \ \ \mbox{on} \ \ \partial K\}\, .
$$
This space is the closure of the space $\mathcal{C}^\infty_c(\overline{Q}\setminus\overline{K})$ with respect to the Dirichlet norm. We also need the two functional spaces of vector fields
$$
\mathcal{V}_{*}(\Omega) = \{ v \in H^1_*(\Omega) \ | \ \nabla \cdot v = 0 \ \ \text{in} \ \Omega \} \ \ \ \text{and} \ \ \ \mathcal{V}(\Omega) = \{ v \in H^{1}_{0}(\Omega) \ | \ \nabla \cdot v = 0 \ \ \text{in} \ \Omega \}.
$$
Assuming \eq{zeroflux2}, we say that a vector field $u \in \mathcal{V}_{*}(\Omega)$ is a \textit{weak solution} of \eqref{nsstokes0}-\eqref{boundaryc0} if $u$
verifies the boundary conditions in the trace sense and
\begin{equation} \label{nstokesdebil}
\eta \int_{\Omega} (\nabla u \cdot \nabla \varphi) \, dx + \int_{\Omega} (u\cdot\nabla)u \cdot \varphi \, dx = 0  \ \ \ \ \ \ \forall \varphi \in \mathcal{V}(\Omega).
\end{equation}

It is well-known \cite[Section IX.4]{galdi2011introduction} that a solution always exists and that it is unique provided that
$\|h\|_{H^{1/2}(\partial Q)}$ is sufficiently small, see also \cite[Section 3]{gazspe} for the particular case of a domain with obstacle.
The flow of the fluid exerts a force $F_{K}$ over the obstacle, which can be computed through the stress tensor $\mathbb{T}(u,p)$
of the fluid \cite[Chapter 2]{landau} and, in a weak sense, is defined as
\begin{equation} \label{weakforce2}
F_{K}(u,p)= - \langle \mathbb{T}(u,p) \cdot \hat{n} , 1 \rangle_{\partial K} \, .
\end{equation}
Here $\langle \cdot , \cdot \rangle_{\partial K}$ denotes the duality pairing between $W^{-\frac{2}{3}, \frac{3}{2}}(\partial K)$ and $W^{\frac{2}{3}, 3}(\partial K)$, while the minus sign is due to the fact that the outward unit normal $\hat{n}$ to $\Omega$ is directed
towards the interior of $K$.
In the case of suspension bridges, the boundary conditions should model an horizontal inflow on the ($2d$ or $3d$) face $x=-L$ of $Q$, as in conditions
\eq{turbulent2d} and \eq{laminar3d}, see also \cite{filippoclara}. Then, the most relevant component of the force \eq{weakforce2}, leading to structural instability,
is the {\em lift force} $\mathcal{L}_{K}(u,p)$ which is oriented vertically and, in our generalized context, can be computed as
$$
\mathcal{L}_{K}(u,p) = F_{K}(u,p) \cdot {\textbf{v}},
$$
where $\textbf{v}$ is the unit vector in the $y$-direction in the 2d space and in the $z$-direction of $z$ in the 3d space.
The connection between the unique solvability of \eqref{nsstokes0}, the existence of symmetric solutions and the appearance of a lift force over $K$ is expressed in the following result:

\begin{proposition}\label{symmetricsolutionns}
Assume \eqref{OmegaR}. For any $h\in H^{1/2}(\partial \Omega)\cap\mathcal{C}(\partial Q)$ satisfying \eqref{zeroflux2}-\eqref{boundaryc0},
there exists a weak solution of \eqref{nsstokes0}. Moreover, there exists $\chi>0$ such that, if
$\|h\|_{H^{1/2}(\partial Q)} <\chi$, then the weak solution is unique. Furthermore, if the obstacle $K$ is symmetric with respect to the $y$-direction (2d case) or to the $z$-direction (3d case) and
the boundary datum satisfies $\|h\|_{H^{1/2}(\partial Q)} <\chi$ and
$$
h(x,-y)=h(x,y) \quad \text{(2d case)} \qquad \text{or} \qquad h(x,y,-z)=h(x,y,z) \quad \text{(3d case)} \quad \text{on} \ \ \partial Q,
$$
then the fluid exerts no lift force on the obstacle, that is, $\mathcal{L}_{K}(u,p)=0$.
\end{proposition}

Note that the symmetry assumption on the inflow is satisfied by a Poiseuille flow but {\em not} by a Couette flow, see \eq{laminar3d}
and \eq{turbulent2d} below. Proposition \ref{symmetricsolutionns} is equally valid if we drop the continuity assumption on $h$; this assumption
is put only for compatibility with the remaining parts of the present work. From \cite{gazspe} we know that the constant $\chi$ in Proposition \ref{symmetricsolutionns} depends on:\par
\noindent
- the viscosity $\eta$, with $\eta\mapsto\chi(\eta)$ being increasing;\par
\noindent
- the geometric measures $L$, $a$, $b$, that modify the embedding constants for $H_{*}^{1}(\Omega),H_{0}^{1}(\Omega)\subset L^{4}(\Omega)$;\par
\noindent	
- the constant $\Gamma$ describing the size of the solenoidal extension, see Theorem \ref{teononconstant}.\par

In \cite{gazspe} the Authors merely considered constant inflows and gave explicit bounds on $\chi$. The main novelty of the present paper
is that we can also handle much more general inflow/outflow problems and still ensure
that no lift force is exerted on the obstacle, as in Proposition \ref{symmetricsolutionns}.

\section{An extension result and two new inflow-outflow models}\label{functional}

We decompose $\partial Q$ as the union of its $(d-1)$-dimensional faces, that we name in dimension $d=2$ by
$$
\Sigma_{1} =\partial Q \cap  \{x=L\} , \ \quad \Sigma_{2} = \partial Q \cap \{y =L\}, \ \quad \Sigma_{3} = \partial Q \cap \{x=-L\}, \ \quad \Sigma_{4} = \partial Q \cap \{y= -L\},
$$
%
%
and in dimension $d= 3$ by
\renewcommand\arraystretch{1}
$$
\begin{array}{lll}
\Sigma_{1} = \partial Q \cap  \{x=L\} , \qquad & \Sigma_{2}  =\partial Q \cap \{y =L\}, \qquad & \Sigma_{3} = \partial Q \cap \{x=-L\}, \\[6pt]
\Sigma_{4} = \partial Q \cap \{y= -L\}, \qquad & \Sigma_{5} = \partial Q \cap  \{z=L\}, \qquad & \Sigma_{6} = \partial Q \cap \{z=-L\}.
\end{array}
$$
We point out that, while in the 2d case we numbered the faces of $\partial Q$ counterclockwise, in the 3d case we kept this ordering and simply added the two extra faces in the $z$-direction. Then, denoting by $h _i 
$ the restriction of $h$ to $\Sigma _i$, the continuity of $h$ at the vertices of $Q$    
in dimension $2$ and at the edges of $Q$ in dimension $3$ 
reads 
\begin{equation}\label{boundaryhcont}
h _ i = h _ j \quad \text{ on } \ \Sigma _i \cap \Sigma _j \,, \qquad \forall i ,j \in \{1,\dots,d\} \,.
\end{equation} 
%
Aim of this section is to construct a vector field $A_1$ as in Step 1 of the outline, namely, a vector field $A_{1} \in H^{1}(Q)\cap \mathcal{C}(\overline{Q})$ such that $A_{1} |_{\partial Q} =h$. We point out that, if the boundary datum is not solenoidal, then the construction of a solenoidal extension to $\Omega_{0}$ cannot be performed merely by the use of
cutoff functions and it is therefore necessary to extend it first to $\overline{Q}$. In view of the special choice of the geometry, such extension can be explicitly found by taking the convex combination of the boundary datum on opposite faces of $\partial Q$.
\newpage
\begin{theorem} \label{gadliardo}
Let $h$  satisfy assumptions \eqref{zeroflux2}.\par\smallskip\noindent	
\begin{itemize}[leftmargin=*]
\item[$\bullet$] For $d= 2$, the function $A_{1} \in H^{1}(Q)\cap \mathcal{C}(\overline{Q})$ defined on  $ \overline{Q}$  by
$$
\begin{aligned}
A_{1}(x,y) & = \dfrac{L+x}{2L} h_{1}(y) + \dfrac{L-x}{2L} h_{3}(y) + \dfrac{L+y}{2L} \left[ h_{2}(x) - \dfrac{L+x}{2L} h_{2}(L) - \dfrac{L-x}{2L} h_{2}(-L) \right] \\[6pt]
& \hspace{4.5mm} + \dfrac{L-y}{2L} \left[ h_{4}(x) - \dfrac{L+x}{2L} h_{4}(L) - \dfrac{L-x}{2L} h_{4}(-L) \right]
\end{aligned}
$$
is an extension of $h$ to $\overline Q$, whose $H ^1$-norm on $Q$ can be explicitly computed in terms of
the $H ^1$-norm of  the functions $h _i$ on $\Sigma _i$, for $i = 1, \dots, 4$.

\item[$\bullet$] For $d = 3$, assuming in addition that $h$ vanishes at the vertices of $Q$, the function
$A_{1}\in H^{1}(Q)\cap \mathcal{C}(\overline{Q})$ defined on $\overline{Q}$ by
$$
\begin{aligned}
A_{1}(x,y) & = \dfrac{L+x}{2L} h_{1}(y,z) + \dfrac{L-x}{2L} h_{3}(y,z) + \dfrac{L+y}{2L} \left[ h_{2}(x,z) - \dfrac{L+x}{2L} h_{2}(L,z) - \dfrac{L-x}{2L} h_{2}(-L,z) \right] \\[6pt]
& \hspace{3.5mm} + \dfrac{L-y}{2L} \left[ h_{4}(x,z) - \dfrac{L+x}{2L} h_{4}(L,z) - \dfrac{L-x}{2L} h_{4}(-L,z) \right] \\[6pt]
& \hspace{3.5mm} + \dfrac{L+z}{2L} \left[ h_{5}(x,y) - \dfrac{L+x}{2L} h_{5}(L,y) - \dfrac{L-x}{2L} h_{5}(-L,y) - \dfrac{L+y}{2L} h_{5}(x,L) - \dfrac{L-y}{2L} h_{5}(x,-L) \right] \\[6pt]
& \hspace{3.5mm} + \dfrac{L-z}{2L} \left[ h_{6}(x,y) - \dfrac{L+x}{2L} h_{6}(L,y) - \dfrac{L-x}{2L} h_{6}(-L,y) - \dfrac{L+y}{2L} h_{6}(x,L) - \dfrac{L-y}{2L} h_{6}(x,-L) \right]
\end{aligned}
$$
is an extension of $h$ to $\overline Q$, whose $H ^1$-norm on $Q$ can be explicitly computed in terms of
the $H ^1$-norm of  the functions $h _i$ on $\Sigma _i$, for $i = 1, \dots, 6$.
\end{itemize}
\end{theorem}
\noindent
\begin{proof}
The result is obtained by computing $A_1$  on $\partial Q$, taking into account
conditions \eqref{boundaryhcont}\,.
\end{proof}

\bigskip
In order to highlight the relevance of Theorem \ref{gadliardo}, we introduce here two new inflow-outflow models  in which
the boundary velocity is not solenoidal. In the 2d case we suggest a model for a turbulent flow,  whereas in the 3d case we suggest a model for an ``almost laminar'' flow.
We emphasize that these models should not be interpreted as a precise description of turbulent or laminar flows in a channel. They  merely serve as possible boundary data to be prescribed in inflow-outflow problems, showing a new level of complexity that can be treated by the methods presented in this article. 

\par

\medskip
{\it A $2d$-model for a turbulent flow.} 
In the planar domain $\Omega_{0}$, we consider an inflow of Couette type and an outflow of ``modified Couette'' type. More precisely, for some $\lambda\gg1$ and $\alpha, \tau > 0$, we take the boundary datum $h : \partial \Omega_{0} \longrightarrow \R ^2$ as
$$
h = 0 \ \text{ on } \Sigma _ 4 \cup \partial K \, ; \qquad
h(x,L)=
\left(\begin{array}{c}
\!\!2\lambda L\!\!\\
0
\end{array}\right)\ \ \forall|x|\le L\, ; \qquad
h(-L,y)=
\left(\begin{array}{c}
\!\!\lambda(y+L)\!\!\\
0
\end{array}\right)\ \ \forall|y|\le L;
$$

\begin{equation} \label{turbulent2d}
h(L,y)=\left(\begin{array}{c}
\lambda(y+L)\\[6pt]
g(y)
\end{array}\right) \ \ \ \forall |y|\le L, \quad \mbox{where} \quad g(y)=\left\{\begin{array}{ll}
0 & \mbox{ if }|y|>b\\
\tau |y|^{\alpha} \sin \left(\dfrac{\pi b}{\sqrt{|y|}} \right) & \mbox{ if }|y|<b\, .
\end{array}\right.
\end{equation}


Conditions \eqref{turbulent2d} aim to model the behavior of the wind on the deck of a bridge, with no-slip condition on the boundary of the bridge and on the floor, see \eqref{turbulent2d}$_1$. Indeed, it is well-known that there is no wind at the ground level (at $y=-L$) and the velocity of the wind increases with altitude, reaching a maximum strength at $y=L$, see \eqref{turbulent2d}$_2$-\eqref{turbulent2d}$_3$. Moreover, this law is linear and, at high altitude, the wind is very strong (large $\lambda$).
The datum is assumed to be regular on the inflow edge $x = - L$, then the flow becomes turbulent after bypassing the obstacle $K$ (see Figure \ref{WT} for the illustration of an experiment), and it 
regularizes only partially  on the outflow edge $x = L$, as in \eqref{turbulent2d}$_4$: in the wake of the obstacle, the flow oscillates as $y\to0$, more quickly but with decreasing amplitude  towards the center of the outflow edge. The intensity of the turbulent motion is measured by the positive parameter $\tau$. 

\begin{figure}[H]
	\begin{center}
		\includegraphics[height=63mm,width=63mm]{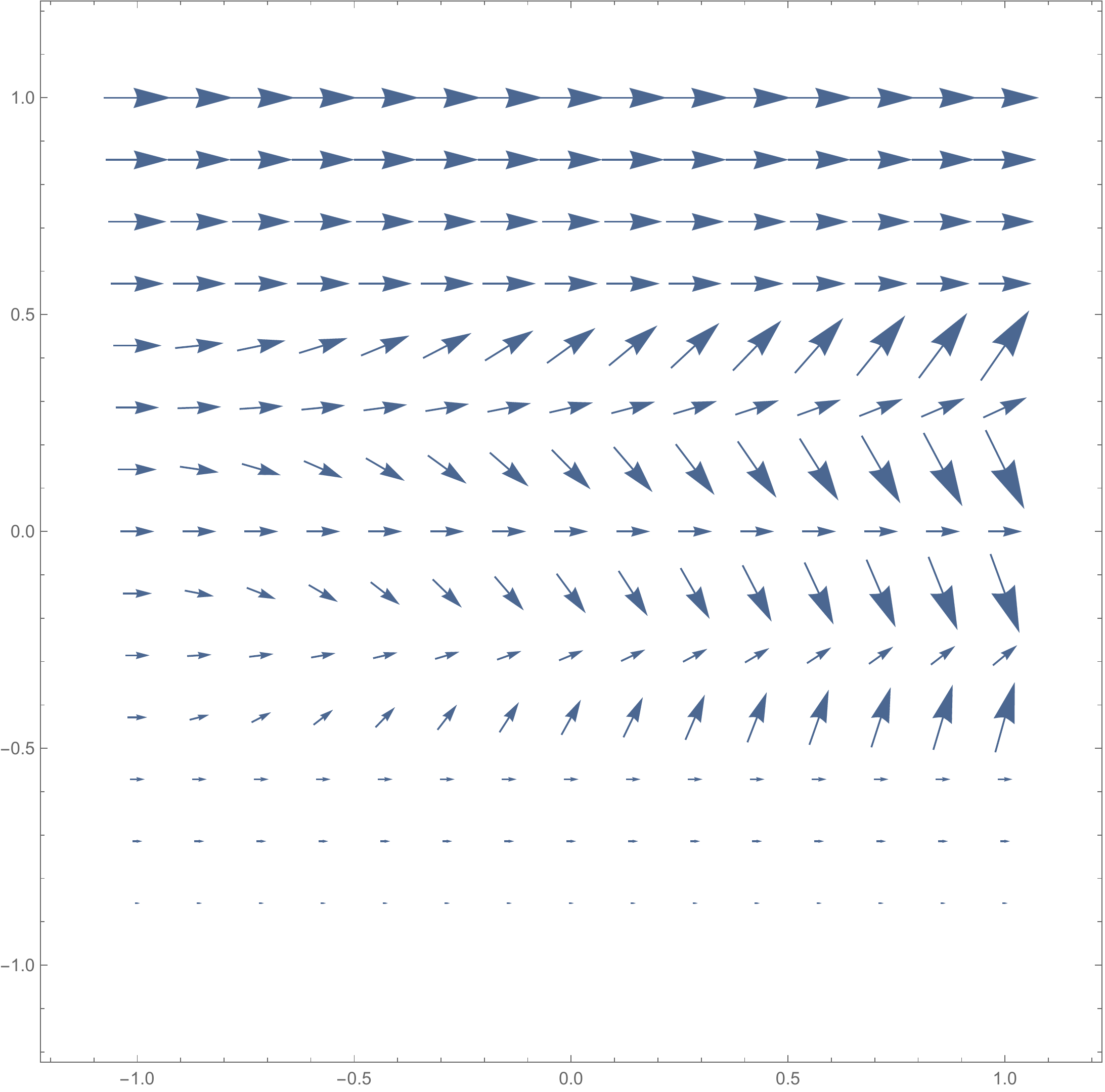}\qquad \qquad \includegraphics[height=63mm,width=63mm]{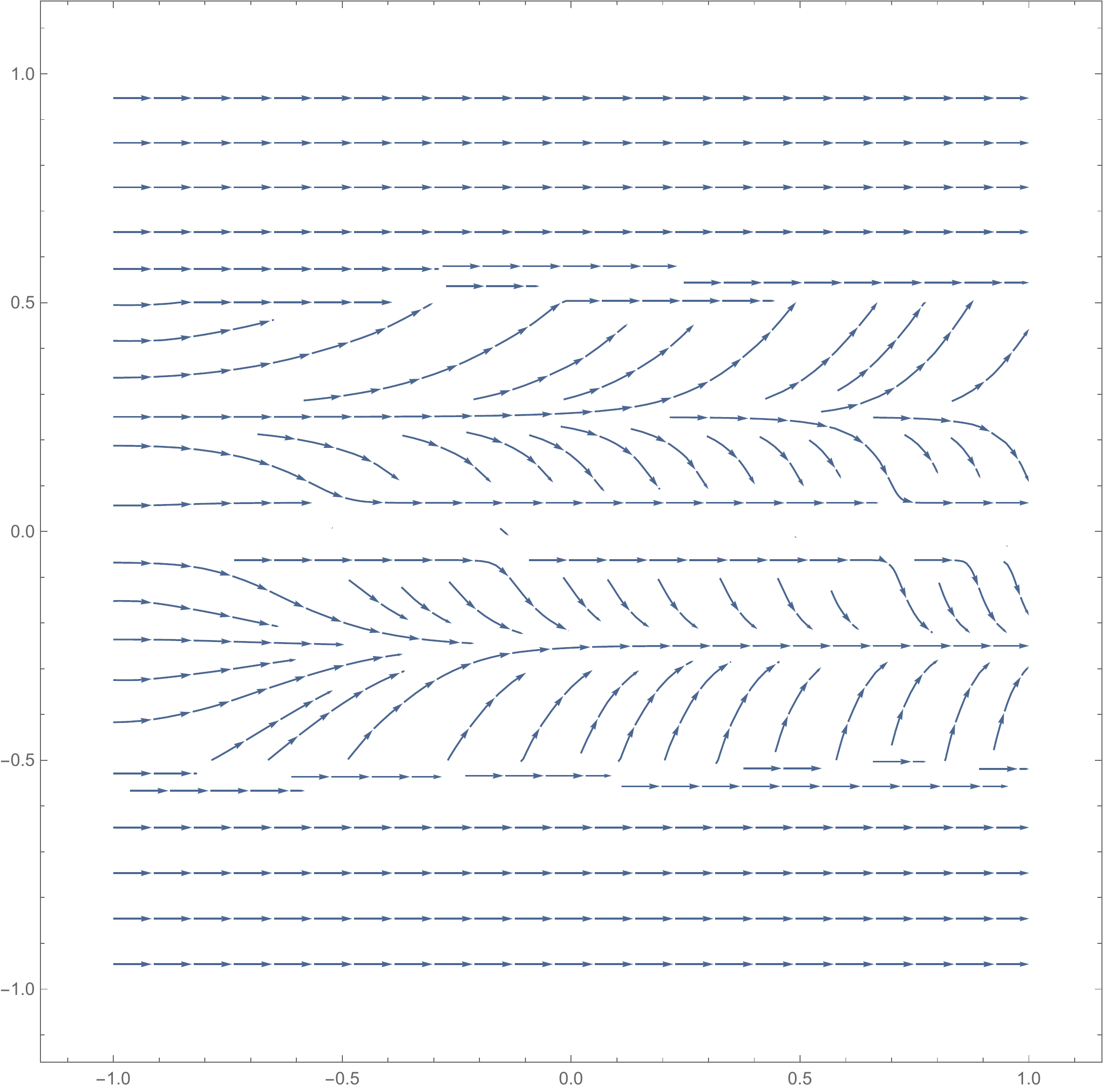}		
	\end{center}
	\vspace*{-3mm}
	\caption{Vector plot (left) and stream plot (right) of the field $A_{1}$ defined in \eqref{estensione2d}.}\label{2dmodel}
\end{figure}

The boundary datum $h$ defined in  \eqref{turbulent2d} satisfies the assumptions in  \eqref{zeroflux2}, and hence it admits a solenoidal extension. The intermediate (non-solenoidal) extension given by Theorem
\ref{gadliardo} reads
\begin{equation} \label{estensione2d}
A_1(x,y)=\left(\begin{array}{c}
\lambda(y+L)\\[6pt]
\dfrac{L+x}{2L}\, g(y)
\end{array}\right)\qquad\forall(x,y)\in \overline{Q}\, .
\end{equation}
The vector and stream plots of the field $A_1$ on $\overline{Q}$ (without the obstacle) are displayed in Figure \ref{2dmodel}, for $L=1$, $a=0.7$, $b=0.5$, $\lambda=3$, $\alpha=0.1$ and $\tau=10$.\par

\medskip
{\it A $3d$-model for an ``almost-laminar'' flow.} 
In the 3d domain $\Omega_{0}$,  we consider an inflow of Poiseuille type and an outflow of ``modified Poiseuille'' type, namely
 we take the datum $h : \partial \Omega_{0} \longrightarrow \R ^3$ as
$$ 
h=0 \ \ \mbox{ on } \ \ \Sigma_{2} \cup \Sigma_{4} \cup  \Sigma_{5} \cup \Sigma_{6} \cup \partial K \, ; \qquad
h(-L,y,z)=\left(\begin{matrix}
2L^2 - y^2 - z^2 \\[3pt]
0\\[3pt]
0
\end{matrix}\right) \qquad \forall (y,z) \in [-L,L]^2\,;
$$
\begin{equation}\label{laminar3d}
h(L,y,z)=\left(\begin{matrix}
2L^2-y^2-z^2\\[8pt]
\dfrac{y(y^2-L^2)(z^2-L^2)}{L^4}\\[8pt]
\dfrac{z(y^2-L^2)(z^2-L^2)}{L^4}
\end{matrix}\right) \qquad \forall (y,z) \in [-L,L]^2\,.
\end{equation}

\begin{figure}[H]
	\begin{center}
	\includegraphics[height=63mm,width=63mm]{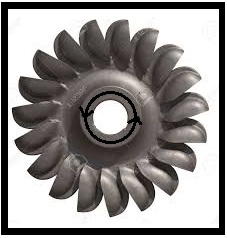}\qquad \qquad \includegraphics[height=63mm,width=63mm]{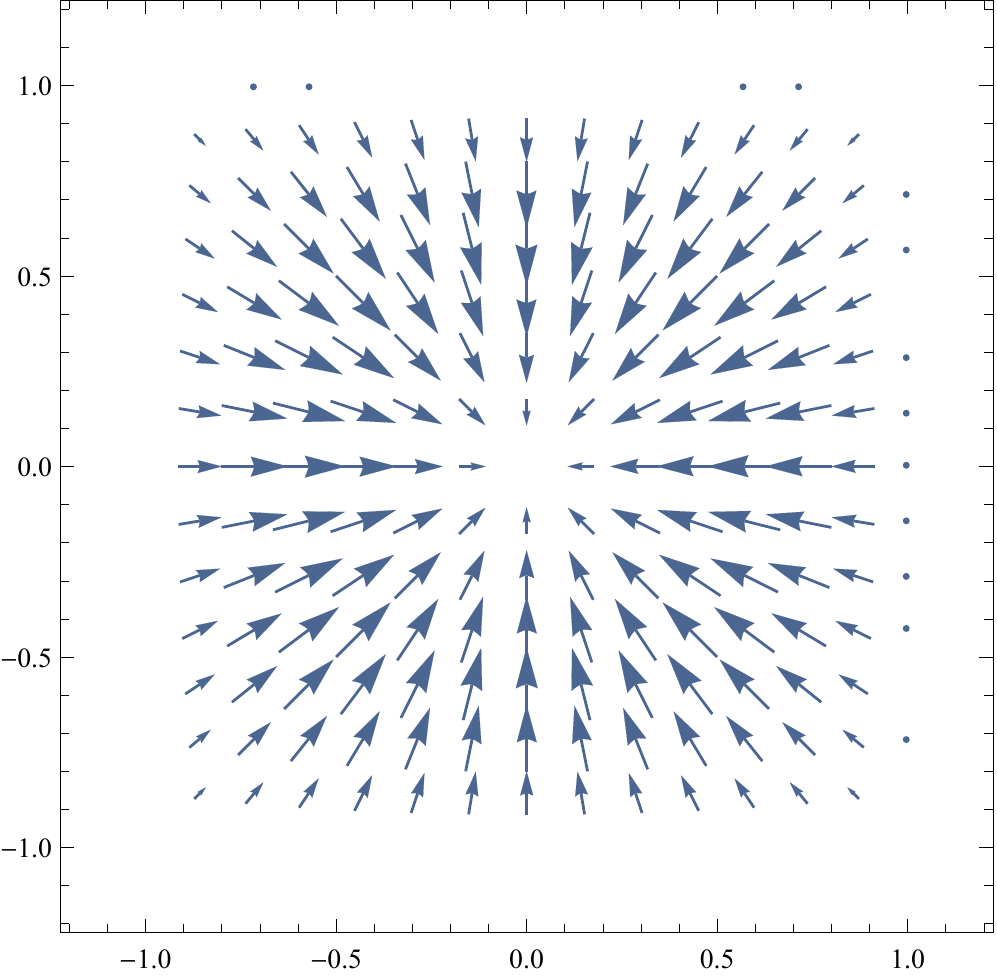}		
	\end{center}
	\vspace*{-3mm}
	\caption{2d view of the inflow-outflow faces of the cube in the 3d model.}\label{3dmodel}
\end{figure}

\begin{figure}[h]
	\begin{center}
	\includegraphics[height=63mm,width=63mm]{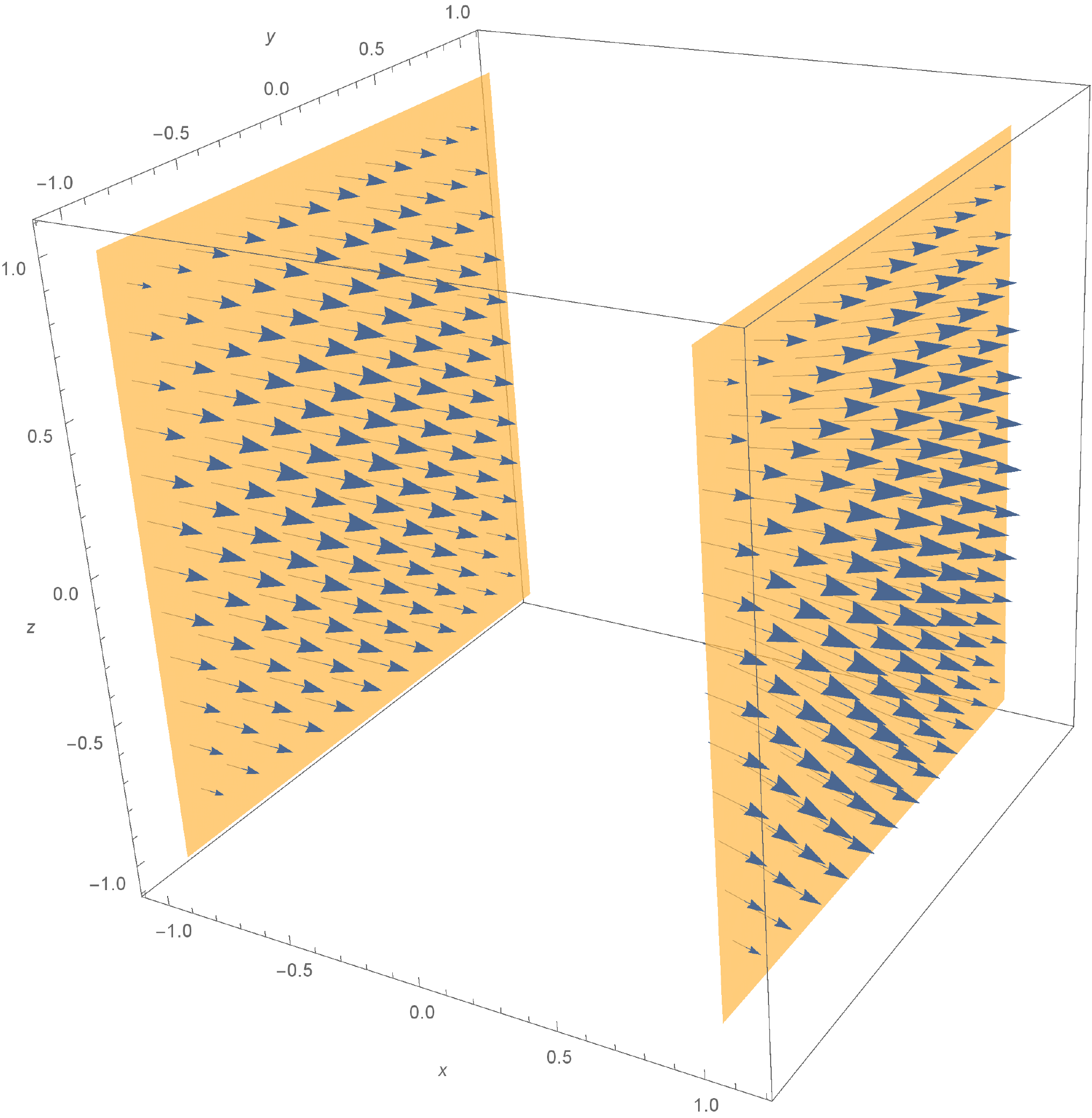}\qquad \qquad \includegraphics[height=63mm,width=63mm]{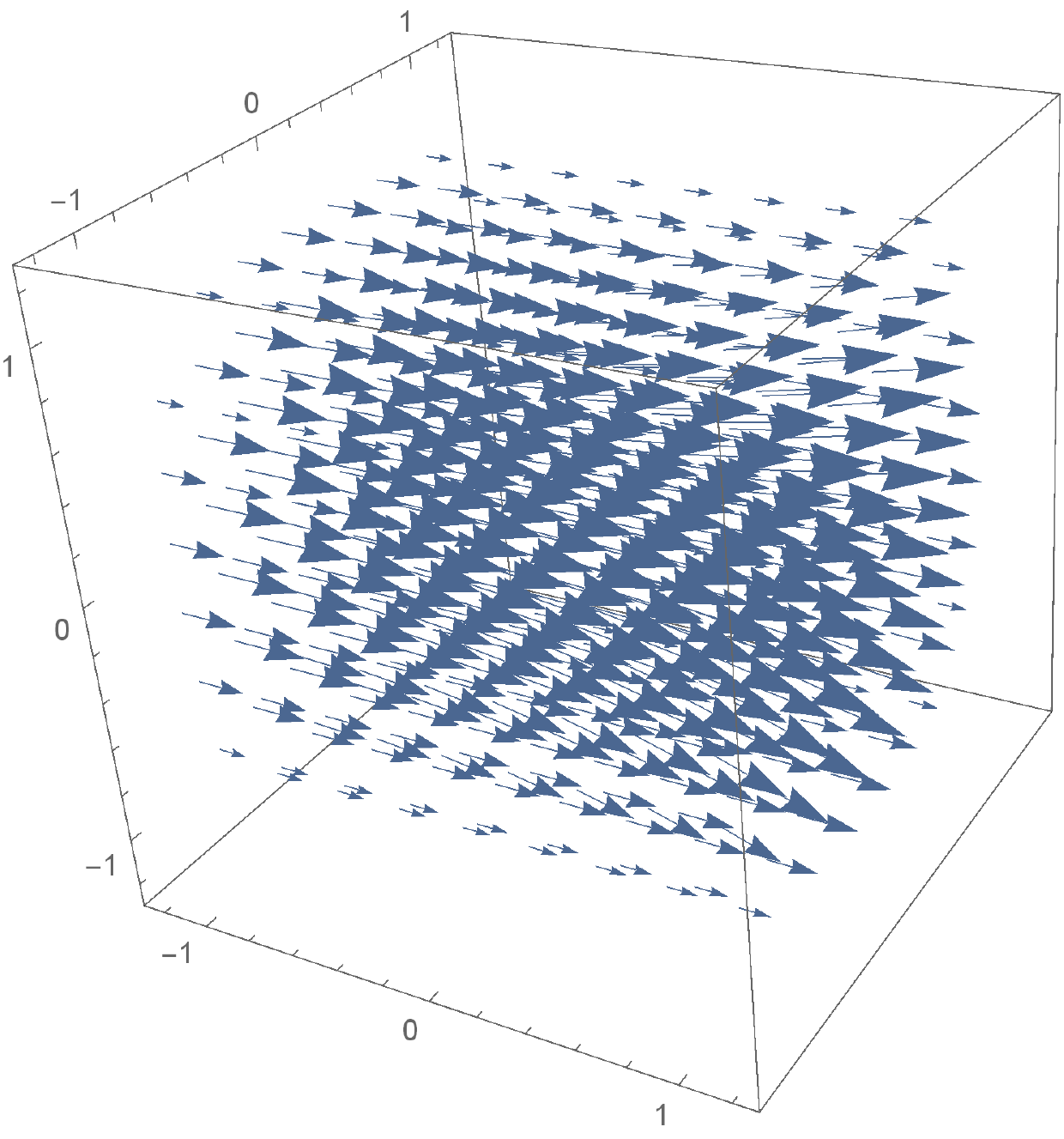}		
	\end{center}
	\vspace*{-3mm}
	\caption{Left: 3d view of the inflow-outflow faces of the cube in the 3d model, for $L=1$. Right: vector plot of the field $A_{1}$ defined in \eqref{estensione3d} for $L=1$.}\label{3dmodelbis}
\end{figure}

Condition \eqref{laminar3d}$_2$ aims at representing a regular inflow, which is expected to slightly modify after by-passing the obstacle, but then tends to recompose at the outflow face, as in \eqref{laminar3d}$_3$ (see \cite{bordogna2020effects} for experimental evidence). In fact, 
in a wind tunnel, the inflow  is generated by a turbine (see the left picture in Figure \ref{3dmodel}), which usually reproduces a Poiseuille flow,
which is faster at the midpoint of the inflow face and vanishes on the edges of this face, as modeled by \eqref{laminar3d}$_2$.
If the inflow is sufficiently small, the motion of the fluid remains almost laminar. It is then reasonable to consider  an outflow which tends to
redistribute regularly also in the wake of the obstacle: this is the reason why in \eqref{laminar3d}$_3$ the flow is oriented
towards the center of the outlet (see Figure \ref{3dmodel} on the right).
The vector field $h$ given by \eqref{laminar3d}$_2$-\eqref{laminar3d}$_3$ is also represented in 3d in the left picture of Figure \ref{3dmodelbis}.
Finally,  on the four remaining faces of the cube (the wind tunnel walls) and on the obstacle, the velocity is zero due to the viscosity:
this is modeled by \eqref{laminar3d}$_1$.

The boundary datum $h$ defined in  \eqref{laminar3d} satisfies the assumptions in  \eqref{zeroflux2}, and hence it admits a solenoidal extension. 
Notice also that $h$ vanishes at the  vertices of the cube $Q$, as requested in \eqref{zeroflux2}. The intermediate (non-solenoidal) extension \eqref{laminar3d}, given by Theorem,
\ref{gadliardo} reads
\begin{equation} \label{estensione3d}
A_1 (x,y,z)=\left(\begin{matrix}
2L^2-y^2-z^2\\[8pt]
\dfrac{x+L}{2L}\, \dfrac{y(y^2-L^2)(z^2-L^2)}{L^4}\\[8pt]
\dfrac{x+L}{2L}\, \dfrac{z(y^2-L^2)(z^2-L^2)}{L^4}
\end{matrix}\right) \qquad \forall (x,y,z) \in \overline{Q}\,.
\end{equation}
The vector plot of the field $A_1$ on $Q$ is displayed in Figure \ref{3dmodelbis} right, for $L=1$.


\section{A Lipschitz function  with gradient of minimal $L^\infty$-norm}\label{sec:phi}

Aim of this section is to construct a scalar function $\phi$ as in Step 2 of the outline.

\begin{itemize}[leftmargin=*]
\item[$\bullet$] For $d = 2$, we denote
by $Q_+$ the intersection of $Q$ with the first quadrant. Setting
$$\gamma (x, y): = (L-a) y + (b-L) x + L ( a-b) \qquad \forall (x,y) \in \mathbb{R}^2,$$
we decompose $Q_+$ as
$Q _+ = \mathcal Q_0 \cup \mathcal Q _1 \cup \mathcal Q _2$,
where
$$
\mathcal Q _0 := Q _+ \cap \overline{P},\quad
\mathcal Q _1 :=  ( Q _+  \setminus \overline{P}) \cap\big \{\gamma \leq 0 \},\quad
\mathcal Q _2 := ( Q _+  \setminus \overline{P})  \cap  \big \{ \gamma \geq 0 \big \},
$$
see Figure \ref{starssezione4} on the left. Then we define $\phi$ on $Q$ as the function given on $Q_+$ by
 \begin{equation}\label{def:phi2}
 \phi (x, y) = \begin {cases}
 0 & \text { if } (x, y) \in \mathcal Q _0
 \\ \noalign{\medskip}
1-\frac{L-x}{L-a} & \text{ if } (x, y) \in \mathcal Q_1
\\ \noalign{\medskip}
1-\frac{L-y}{L-b}  & \text { if } (x, y) \in \mathcal Q _2\,.
\end{cases}
\end{equation}
and extended by even reflection to the other quadrants.  We have
\begin{equation}\label{f:phi2}
\| \phi \| _{L ^ \infty (Q)}  = 1 \qquad \text{ and } \qquad \| \nabla \phi \| _{L ^ \infty (Q)}  =  \max \Big \{  \frac{1}{L-a}, \frac{1}{L-b}  \Big \}
= \frac{1}{L-a}
\,.
\end{equation}

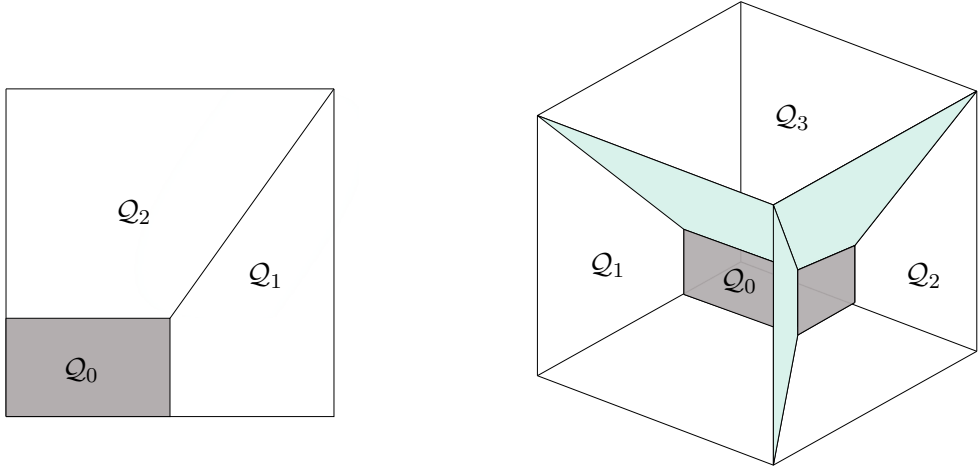
\begin{figure}[H]
	\begin{minipage}{0.5\linewidth}
		\vskip - 3 cm
		\centering
		\def\svgwidth{7cm}
		\begingroup%
		\makeatletter%
		\providecommand\color[2][]{%
			\errmessage{(Inkscape) Color is used for the text in Inkscape, but the package 'color.sty' is not loaded}%
			\renewcommand\color[2][]{}%
		}%
		\providecommand\transparent[1]{%
			\errmessage{(Inkscape) Transparency is used (non-zero) for the text in Inkscape, but the package 'transparent.sty' is not loaded}%
			\renewcommand\transparent[1]{}%
		}%
		\providecommand\rotatebox[2]{#2}%
		\ifx\svgwidth\undefined%
		\setlength{\unitlength}{377.58676bp}%
		\ifx\svgscale\undefined%
		\relax%
		\else%
		\setlength{\unitlength}{\unitlength * \real{\svgscale}}%
		\fi%
		\else%
		\setlength{\unitlength}{\svgwidth}%
		\fi%
		\global\let\svgwidth\undefined%
		\global\let\svgscale\undefined%
		\makeatother%
		\begin{picture}(1,1.38846134)%
		\put(0,0)
		{\includegraphics[height=55mm]{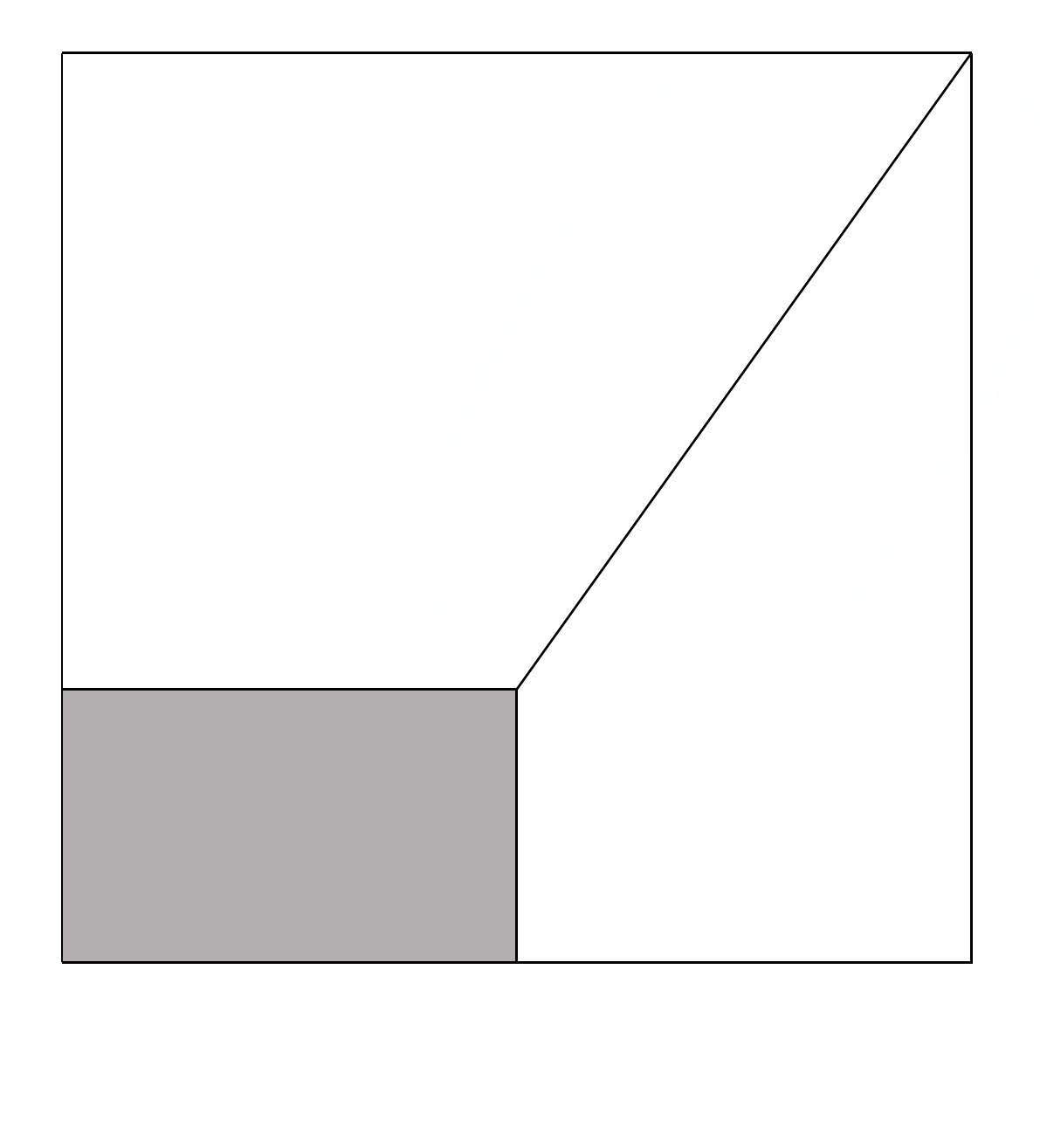}
			\hskip 2 cm
			\includegraphics[height=65mm]{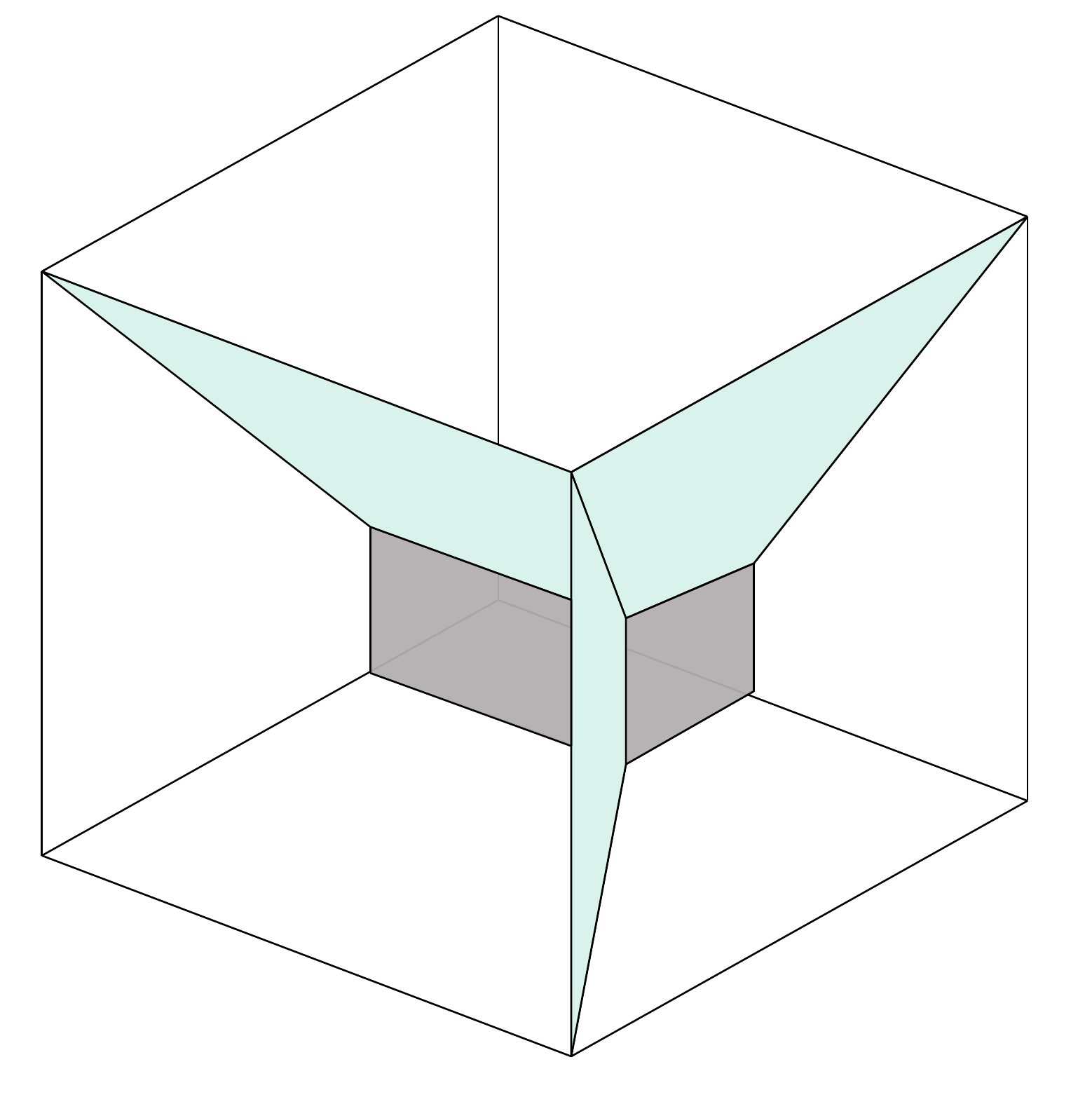}}
		\put(0.15,0.20){\color[rgb]{0,0,0}\makebox(0,0)[lb]{\smash{$\mathcal Q_0$}}}%
		\put(0.25,0.50){\color[rgb]{0,0,0}\makebox(0,0)[lb]{\smash{$\mathcal Q_2$}}}%
		\put(0.5,0.38){\color[rgb]{0,0,0}\makebox(0,0)[lb]{\smash{$\mathcal Q_1$}}}%
		\put(1.4,0.37){\color[rgb]{0,0,0}\makebox(0,0)[lb]{\smash{$\mathcal Q_0$}}}%
		\put(1.15,0.40){\color[rgb]{0,0,0}\makebox(0,0)[lb]{\smash{$\mathcal Q_1$}}}%
		\put(1.75,0.38){\color[rgb]{0,0,0}\makebox(0,0)[lb]{\smash{$\mathcal Q_2$}}}%
		\put(1.5,0.68){\color[rgb]{0,0,0}\makebox(0,0)[lb]{\smash{$\mathcal Q_3$}}}%
		
		\end{picture}%
		\endgroup%
	\end{minipage}%
	\begin{minipage}{0.5\linewidth}
		\centering
		\def\svgwidth{7cm}
		\begingroup%
		\makeatletter%
		\providecommand\color[2][]{%
			\errmessage{(Inkscape) Color is used for the text in Inkscape, but the package 'color.sty' is not loaded}%
			\renewcommand\color[2][]{}%
		}%
		\providecommand\transparent[1]{%
			\errmessage{(Inkscape) Transparency is used (non-zero) for the text in Inkscape, but the package 'transparent.sty' is not loaded}%
			\renewcommand\transparent[1]{}%
		}%
		\providecommand\rotatebox[2]{#2}%
		\ifx\svgwidth\undefined%
		\setlength{\unitlength}{379.39037574bp}%
		\ifx\svgscale\undefined%
		\relax%
		\else%
		\setlength{\unitlength}{\unitlength * \real{\svgscale}}%
		\fi%
		\else%
		\setlength{\unitlength}{\svgwidth}%
		\fi%
		\global\let\svgwidth\undefined%
		\global\let\svgscale\undefined%
		\makeatother%
		\endgroup%
	\end{minipage}
	\caption{Decomposition of $Q_+$ into the regions $\mathcal Q_i$, in $2d$ (left) and in $3d$ (right).}
	\label{starssezione4}
\end{figure}

\item[$\bullet$]  For $d = 3$, we denote
by $Q_+$ the intersection of $Q$ with the first octant. Setting
$$
\begin{array}{cc}
\displaystyle  \gamma _ 1 ( y, z) := (L-c) y + ( b-L) z + L (c-b), \qquad \displaystyle \gamma _ 2 (x, z) := (L-a) z + (c-L) x + L ( a-c), 
\\  \noalign{\medskip}
\displaystyle  \gamma_ 3 (x, y) := (L-b) x + (a-L) y + L (b-a) \qquad \forall (x,y,z) \in \mathbb{R}^3 \, ,
\end{array}
$$
we decompose $Q_+$ as
$Q _+ = \mathcal Q_0 \cup \mathcal Q _1 \cup \mathcal Q _2 \cup \mathcal Q _ 3$,
where (see Figure \ref{starssezione4} on the right)
$$
\begin{array}{ll}
\displaystyle \mathcal Q _0 := Q _+  \cap \overline P
\qquad
& \displaystyle \mathcal Q _1 := (Q _+  \setminus \overline{P}) \cap  \big \{ \gamma_2  \leq 0 \big \} \cap \big  \{ \gamma_3  \geq 0 \big  \}
\\  \noalign{\medskip}
\displaystyle \mathcal Q _2 := (Q _+  \setminus \overline{P}) \cap\big \{\gamma_1 \geq 0 \big \} \cap \big  \{ \gamma_3 \leq 0 \big  \}
\qquad
& \displaystyle \mathcal Q _3 := (Q _+  \setminus \overline{P})  \cap \big \{ \gamma_2  \geq 0 \big \} \cap \big  \{ \gamma_1 \leq 0  \big  \} \, .
\end{array}
$$
Then we define $\phi$ on $Q$ as the function given on $Q_+$ by
\begin{equation}\label{def:phi3}
\phi (x, y, z) := \begin {cases}
0 & \text{ if } (x, y, z) \in \mathcal Q _0
 \\ \noalign{\medskip}
1- \frac{L-x}{L-a} & \text{ if } (x, y, z) \in \mathcal Q _1
\\ \noalign{\medskip}
1 - \frac{L-y}{L-b}  & \text { if }  (x, y, z) \in \mathcal Q _2
\\ \noalign{\medskip}
1 - \frac{L-z}{L-c}  & \text { if }  (x, y, z) \in \mathcal Q _3
\end{cases}
\end{equation}
and extended by even reflection to the other octants.   We have
\begin{equation}\label{f:phi3}
\| \phi \| _{L ^ \infty (Q)}  = 1 \qquad \text{ and } \qquad \| \nabla \phi \| _{L ^ \infty (Q)}  =  \max \left\{  \frac{1}{L-a}, \frac{1}{L-b} , \frac{1}{L-c}  \right\}
= \frac{1}{L-a}
\,.
\end{equation}
\end{itemize}
%


The next result shows that the norm of the function $\phi$  in $W ^ {1, \infty} (Q)$ is
minimal in the class of functions in $W ^ {1, \infty} (Q)$ which are equal to $1$ on $\partial Q$ and vanish on $\overline P$.

\begin{theorem}\label{thm:infty}
The function $\phi$ defined for $d= 2$ by \eqref{def:phi2} and for $d=3$ by \eqref{def:phi3} solves the variational problem
$$\min \Big \{ \| \nabla u \| _{ L ^\infty (\Omega_{0} ) }  \ : \ u \in W ^ {1, \infty}  ( \Omega_{0} ) \,, \ u = 1 \ \text{ on } \ \partial Q  \, ,\  u = 0 \ \text{ on  } \ \partial P \Big \} \,.$$
\end{theorem}\noindent
\proof  Let $w$ be the
infinity-harmonic potential of
$\overline P$ relative to $Q$, namely  the unique viscosity solution to
the boundary value problem
\begin{equation}\label{f:pot}
\begin{cases}
-\Delta_\infty  w = 0, & \text{in}\ \Omega_{0},\\
w  = 1, &\text{on}\ \partial Q,\\
w  = 0, &\text{on}\  \partial P \, ,
\end{cases}
\end{equation}
where $\Delta _\infty$ denotes the infinity-Laplacian operator, defined for smooth functions $u$ by
$$
\Delta _\infty u\doteq D ^ 2 u \cdot \nabla u \cdot \nabla u \, .
$$
The existence and uniqueness of a viscosity solution to problem \eqref{f:pot} is due to Jensen (see \cite[Section 3]{jensen}), who also proved that $w$ has the following variational property (usually referred to as AML, i.e.\  {\it absolutely minimizing Lipschitz extension}): for every open bounded set $A \subset \Omega_{0} $, and for every function $v \in \mathcal{C} (\overline A)$ such that $ v = w $ on $\partial A$, it holds $\| \nabla w  \| _{L ^\infty (A)} \leq \| \nabla v \|  _{ L ^\infty(A)}$ (see  also
\cite{CEG}).
Then the statement of the theorem is equivalent to
assert that
\begin{equation}\label{f:opt2}
\|\nabla  w  \| _ {L ^\infty (\Omega_{0} ) }  =   \|\nabla  \phi  \| _ {L ^\infty (\Omega_{0} ) }   \,.
\end{equation}
Since $\phi = w$ on $\partial \Omega_{0} $, the inequality $\|\nabla  w  \| _ {L ^\infty (\Omega_{0} ) }  \leq   \|\nabla  \phi  \| _ {L ^\infty (\Omega_{0} ) } $ follows directly from the AML property of $w$. To prove the converse inequality, we recall from \cite[Proposition 9]{CF} that,
if $q\in\partial Q $ and $p\in \partial P $ are two points
such that $|q-p| = {\rm dist}(\overline P, \partial Q)$, then $w$ is affine on the segment $[q,p]$, and hence it agrees with the function $\phi$ defined in  
\eqref{def:phi2}-\eqref{def:phi3}. Therefore,
$$\begin{array}{ll}
\|\nabla w \|_ {L ^\infty (\Omega_{0} )}  &\displaystyle \geq \sup  \left\{ \frac{|w (q) - w (p)| }{|q-p| }  \ :\ q \in \partial Q , \ p \in \partial P\, , \ |q-p| = {\rm dist}(\overline P, \partial Q) \right\}  
\\  \noalign{\medskip}
&\displaystyle =  \sup  \left\{ \frac{|\phi (q) - \phi (p)| }{|q-p| }  \ :\ q \in \partial Q , \ p \in \partial P\, , \ |q-p| = {\rm dist}(\overline P, \partial Q) \right\}  
=   \|\nabla  \phi  \| _ {L ^\infty (\Omega_{0} ) } \,.
\end{array}
$$
This shows that $\phi$ solves the variational problem and completes the proof.
\qed

\section{An upper bound for the Bogovskii constant} \label{duransec}

Aim of this section is to construct a vector field $A_3$ as in Step 3 of the outline.   To that purpose, let us recall that 
the \textit{Bogovskii problem} in $\Omega_{0}$ consists in finding a positive constant $C$,  depending only on $\Omega_{0}$, such that	
\begin{equation}\label{f:bp}
	\forall g \in L ^ 2 _0 (\Omega_{0})\, , \ \exists v \in  H_{0}^{1}(\Omega_{0}) \quad :\quad \nabla\cdot v=g \  \ \text{in} \ \ \Omega_{0} \quad \text{ and } \quad \|\nabla v\|_{L^2(\Omega_{0})} \leq C  \|g\|_{L^2(\Omega_{0})}\, .
	\end{equation}
The smallest among such positive constants $C$, 
that we denote by $C_B (\Omega_{0})$, is the \textit{Bogovskii constant} of $\Omega_{0}$,
see \eqref{bogo}. We now provide an explicit upper bound for $C _B (\Omega_{0})$.
Notice that proving an inequality of the form $C _B (\Omega_{0}) \leq M$ is equivalent to proving that the claim in \eqref{f:bp} holds true with $C = M$. Let us introduce the following notation:
\begin{itemize}[leftmargin=*]
\item{}
for $d= 2$, we set
{\footnotesize
\begin{equation}\label{sigmagamma}
\begin{aligned}
& \hspace{-3.5mm} \sigma_{2}  =  \dfrac{(L-a)^2}{16 a(L+a)} \left[ (L+3a)(L+a-2b) - (L-a) \sqrt{(L+a-2b)^{2} + 8a(L+a)} \right] +3 L^2  - (a + b) L  -a b , \\[6pt]
& \hspace{-3.5mm} \gamma_{2} = \dfrac{(L-a)^2}{8a(L+a)} \left[ (L+3a)(L+a-2b) - (L-a) \sqrt{(L+a-2b)^{2} + 8a(L+a)} \right] + 2L ^ 2 - 2(a+b) L + 2ab  \,;
\end{aligned}
\end{equation}}
\item{} for $d = 3$, we set
\begin{equation}\label{sigmagamma3d}
\begin{aligned}
& \hspace{-2.5mm} \sigma_{3}  =    7 L^3 -  (a + b + c) L^2 - (a b  +a c + b c) L -a b c , \\[6pt]
& \hspace{-2.5mm} \gamma_{3} =   6 L^3- 2 (a + b + c) L^2 - 2 (a b + a c + b c) L+6 a b c \,.
\end{aligned}
\end{equation}
\end{itemize}
We point out that $\sigma_{d}= |\Omega _1| = |\Omega_2| $ and $\gamma_{d} = |\Omega _1 \cap \Omega_2|$, where $\Omega_1$ and $\Omega_2$ are the two domains in which we subdivide $\Omega_{0}$ in order to obtain Theorem \ref{unionstarshapedcor} below, as outlined in Step 3 of
Section \ref{outline}. Note also that, while \eqref{sigmagamma3d} is symmetric in $a, b, c$, this is not the case for 
\eqref{sigmagamma}:  the reason is the different decomposition we performed in $2d$ and $3d$  (see Figure \ref{stars}). 
This will become  fully clear after reading the proof given hereafter.
\begin{theorem} \label{unionstarshapedcor} There holds $C _B (\Omega_{0}) \leq M$, where  the explicit value of the constant $M$ is given below:\par
\noindent
$\bullet$  for $d=2$, 
letting $\sigma_{2}$ and $\gamma_{2}$ be defined in \eqref{sigmagamma}, 
$$ M \doteq 2  \sqrt{ 2  \left( 1 + \dfrac{8}{\gamma_{2}} (L^2 - ab) \right) } \Bigg[ 129.35 + \dfrac{143.86 \sqrt{\sigma_{2}}}{L-a} + \dfrac{45.36\sigma_{2}}{(L-a)^2}
+ \dfrac{64L^{2}}{(L-a)^2} \left( 13.79 + \dfrac{7.28 \sqrt{\sigma_{2}}}{L-a} \right)^{2} \Bigg]^{1/2} \, ;
$$\par
\noindent
$\bullet$ for $d = 3$, letting $\sigma_{3}$ and $\gamma_{3}$ be defined in \eqref{sigmagamma3d},
{\small
$$
M \doteq \sqrt{ 12  \left( 1 + \dfrac{16}{\gamma_{3}} (L^3 - abc) \right) } \Bigg[ 327.23 + \dfrac{445.17 \sqrt{\sigma_{3}}}{(L-a)^{3/2}} + \dfrac{153.85 \sigma_{3}}{(L-a)^3}
+ \dfrac{144 L^{2}}{(L-a)^2} \left( 22.4 + \dfrac{15.79 \sqrt{\sigma_{3}}}{(L-a)^{3/2}}  \right)^{2} \Bigg]^{1/2} .
$$}
\end{theorem}

In order to prove Theorem \ref{unionstarshapedcor},  we need as a preliminary result the following estimate for the Bogovskii constant of a domain which is star-shaped with respect to a ball.

\begin{proposition}\label{p:starshaped}
Let  $\mathcal{O}\subset\mathbb{R}^d$ be a bounded domain which is star-shaped with respect to a ball $\mathcal{B}\subset\mathcal{O}$
of radius $r>0$.  Then the  following upper bound for the Bogovskii constant $C_B (\mathcal O)$ holds:
$$
\left\lbrace
\begin{aligned}
& C_B (\mathcal O) \leq 2 \Bigg[ 129.35 + \dfrac{71.93}{r} \sqrt{| \mathcal{O} |} + \dfrac{11.34}{r^2} | \mathcal{O} |  + 2 \dfrac{\delta(\mathcal{O})^2}{r^2}\left( 13.79 + \dfrac{3.64}{r} \sqrt{| \mathcal{O} |} \right)^{2} \Bigg]^{1/2} \quad \text{if} \ \ \ d=2, \\[6pt]
& C_B (\mathcal O) \leq \sqrt{6} \Bigg[ 327.23 + \dfrac{157.92}{r^{3/2}} \sqrt{| \mathcal{O} |} + \dfrac{19.23}{r^3} | \mathcal{O} |  + 3 \dfrac{\delta(\mathcal{O})^2}{r^2}\left( 22.4 + \dfrac{5.58}{r^{3/2}} \sqrt{| \mathcal{O} |} \right)^{2} \Bigg]^{1/2} \quad \text{if} \ \ \ d=3,
\end{aligned}
\right.
$$
where $|\mathcal O|$ and  $\delta(\mathcal{O})$ denote, respectively, the Lebesgue measure and the diameter of $\mathcal{O}$.\par
\end{proposition}
\noindent
\proof
After a translation we may assume that the ball $\mathcal{B}$ is centered at the origin of $\mathbb{R}^d$. Let $\omega_{d} \in \mathcal{C}^{\infty}_{0}(\mathbb{R}^{d})$ be the standard radial mollifier whose support coincides with $\mathcal{B}$, that is,
\neweq{mollifier}
\omega_{2}(x,y) =
\left\lbrace
\begin{aligned}
& \dfrac{\ell_2}{r^{2}} \exp\left( \dfrac{r^2}{x^2 + y^2 - r^2} \right) \ \ \ \ & \text{if} \ \ x^2 + y^2 < r^2, \\[6pt]
& 0 & \text{if} \ \ x^2 + y^2 \geq r^2,
\end{aligned}
\right.
\endeq
and
\neweq{mollifier3d}
\omega_{3}(x,y,z) =
\left\lbrace
\begin{aligned}
	& \dfrac{\ell_3}{r^{3}} \exp\left( \dfrac{r^2}{x^2 + y^2 + z^2 - r^2} \right) \ \ \ \ & \text{if} \ \ x^2 + y^2 + z^2 < r^2, \\[6pt]
	& 0 & \text{if} \ \ x^2 + y^2 + z^2 \geq r^2.
\end{aligned}
\right.
\endeq
In \eqref{mollifier}-\eqref{mollifier3d}, $\ell_{d} > 0$ is the normalization constant such that $\|\omega_{d}\|_{L^1(\mathcal{B})} = 1$; hence,
\begin{equation} \label{normalcons}
\ell_{2}=\left(2\pi\int_0^1 t \, e^{1/(t^2-1)}dt\right)^{-1}\approx 2.14357, \qquad \ell_{3}=\left(4\pi\int_0^1 t^{2} e^{1/(t^2-1)}dt\right)^{-1}\approx 2.26712.
\end{equation}
Given $g \in L^{2}_0(\mathcal{O})$, Bogovskii \cite{bogovskii1979solution} showed that a solution $W \in H_{0}^{1}(\mathcal{O})$ of the problem $\nabla \cdot W = g$ can be written as
\begin{equation} \label{bog1}
W(\xi) = \int\limits_{\mathcal{O}} \int\limits_{0}^{1} \dfrac{\xi - \xi'}{t^{3}}\, \omega_{d} \left( \xi' + \dfrac{\xi - \xi'}{t} \right) g(\xi') \ dt\,d\xi' \ \ \ \ \ \forall \xi \in \mathbb{R}^{d}.
\end{equation}
Following \cite{duran2012elementary}, we differentiate \eqref{bog1} under the integral sign, and we interpret the partial derivatives of the field
$W= ( W ^ 1, \dots, W ^ d)$ as operators acting on the function $g$. In other words, for $k,j \in \{1,\ldots,d\}$, 

\neweq{partialvk}
\dfrac{\partial W^{k}}{\partial \xi_{j}}(\xi) = T_{kj,1}(g)(\xi) - T_{kj,2}(g)(\xi)\ \ \ \ \ \forall \xi \in \mathbb{R}^{d}\, , \endeq
where, if $g$ is extended by zero outside $\mathcal{O}$,
$$
\left\lbrace
\begin{aligned}
& T_{kj,1}(g)(\xi) = \lim\limits_{\varepsilon \to 0} \int\limits_{\varepsilon}^{1} \int\limits_{\mathbb{R}^{d}} \dfrac{1}{t^{2}}  \dfrac{\partial}{\partial \xi_{j}} \left[ \left( \xi'_{k} + \dfrac{\xi_{k} - \xi'_{k}}{t} \right) \omega_{d} \left( \xi' + \dfrac{\xi - \xi'}{t} \right) \right] g(\xi') \ dt\,d\xi', \\[6pt]
& T_{kj,2}(g)(\xi) = \lim\limits_{\varepsilon \to 0} \int\limits_{\varepsilon}^{1} \int\limits_{\mathbb{R}^{d}} \dfrac{\xi'_{k}}{t^{2}} \dfrac{\partial}{\partial \xi_{j}} \left[ \omega_{d} \left( \xi' + \dfrac{\xi - \xi'}{t} \right) \right] g(\xi') \ dt\,d\xi'\,.
\end{aligned}
\right.
$$
In view of \eq{partialvk}, by applying Young's  inequality, we get
\begin{equation}\label{f:Young} 
\| \nabla W \|^{2}_{L^{2}(\mathcal{O})} \leq 2 \left( \sum_{k,j=1}^{d} \left\| T_{kj,1}(g) \right\|^{2}_{L^{2}(\mathcal{O})} + \sum_{k,j=1}^{d} \left\| T_{kj,2}(g) \right\|^{2}_{L^{2}(\mathcal{O})} \right).
\end{equation}
In order to estimate the right hand side of \eqref{f:Young} we recall from \cite[Theorem 3.1]{duran2012elementary} that, for $k,j \in \{1,...,d\}$,
\begin{equation} \label{estduran1}
\left\lbrace
\begin{aligned}
& \| T_{kj,1}(g) \|_{L^{2}(\mathcal{O})} \leq \left( 2^{\frac{d-1}{2}} A_{kj} + 2^{\frac{d}{2}} \widetilde{A}_{kj} \sqrt{ | \mathcal{O} |} \right) \| g \|_{L^{2}(\mathcal{O})}, \\[6pt]
& \| T_{kj,2}(g) \|_{L^{2}(\mathcal{O})} \leq \delta(\mathcal{O}) \left( 2^{\frac{d-1}{2}} B_{kj} + 2^{\frac{d}{2}} \widetilde{B}_{kj} \sqrt{ | \mathcal{O} |} \right) \| g \|_{L^{2}(\mathcal{O})} \,,
\end{aligned}
\right.
\end{equation}
where the constants $A_{kj}$, $\widetilde{A}_{kj}$, $B_{kj}$ and $\widetilde{B}_{kj}$ are explicitly given by
\begin{equation} \label{estduran2}
\left\lbrace
\begin{aligned}
& A_{kj} = \dfrac{1}{r} \| \xi_{k} \, \omega_{d} \|_{L^{1}(\mathcal{B})} + r \left\| \dfrac{\partial^{2}}{\partial \xi^{2}_{j}} (\xi_{k} \, \omega_{d}) \right\|_{L^{1}(\mathcal{B})} \ \ \ \ \widetilde{A}_{kj} = \left\| \dfrac{\partial}{\partial \xi_{j}} (\xi_{k} \, \omega_{d}) \right\|^{\frac{1}{2}}_{L^{1}(\mathcal{B})} \left\| \dfrac{\partial}{\partial \xi_{j}} (\xi_{k} \, \omega_{d}) \right\|^{\frac{1}{2}}_ {L^{\infty}(\mathcal{B})}, \\[6pt]
& B_{kj} = \dfrac{1}{r} \| \omega_{d} \|_{L^{1}(\mathcal{B})} + r \left\| \dfrac{\partial^{2} \omega_{d}}{\partial \xi^{2}_{j}} \right\|_{L^{1}(\mathcal{B})} \ \ \ \ \ \ \ \widetilde{B}_{kj} = \left\| \dfrac{\partial \omega_{d}}{\partial \xi_{j}} \right\|^{\frac{1}{2}}_{L^{1}(\mathcal{B})} \left\| \dfrac{\partial \omega_{d}}{\partial \xi_{j}} \right\|^{\frac{1}{2}}_ {L^{\infty}(\mathcal{B})}\,.
\end{aligned}
\right.
\end{equation}
At this point we  distinguish between the cases $d=2$ and $d=3$:

\medskip
\noindent
$\bullet$ For $d=2$, the constants in \eqref{estduran2} admit the following upper bounds (see Section \ref{apA}):
$$
\begin{aligned}
& A_{11} = A_{22} < 7.29, \ \ \ \ \ \ A_{12} = A_{21} < 3.39, \ \ \ \ \ \ \widetilde{A}_{11} = \widetilde{A}_{22} < \dfrac{1.19}{r}, \ \ \ \ \ \ \widetilde{A}_{12} = \widetilde{A}_{21} < \dfrac{0.66}{r},\\[3pt]
& B_{11} = B_{12} = B_{21} = B_{22} < \dfrac{9.75}{r}, \ \ \ \ \ \ \ \widetilde{B}_{11} = \widetilde{B}_{12} = \widetilde{B}_{21} = \widetilde{B}_{22} < \dfrac{1.82}{r^{2}}.
\end{aligned}
$$
By inserting these values into \eqref{estduran1} we obtain:
\begin{equation} \label{estduran3}
\left\lbrace
\begin{aligned}
& \| T_{11,1}(g) \|_{L^{2}(\mathcal{O})} = \| T_{22,1}(g) \|_{L^{2}(\mathcal{O})} \leq \left( 10.31 + \dfrac{2.38}{r} \sqrt{| \mathcal{O} |} \right) \| g \|_{L^{2}(\mathcal{O})}, \\[6pt]
& \| T_{12,1}(g) \|_{L^{2}(\mathcal{O})} = \| T_{21,1}(g) \|_{L^{2}(\mathcal{O})} \leq \left( 4.8 + \dfrac{2.38}{r} \sqrt{| \mathcal{O} |} \right) \| g \|_{L^{2}(\mathcal{O})}, \\[6pt]
& \| T_{kj,2}(g) \|_{L^{2}(\mathcal{O})} \leq \dfrac{\delta(\mathcal{O})}{r} \left( 13.79 + \dfrac{3.64}{r} \sqrt{| \mathcal{O} |} \right) \| g \|_{L^{2}(\mathcal{O})} \ \ \ \ \ \ \text{for every} \ k,j \in \{1,2\}.
\end{aligned}
\right.
\end{equation}
\noindent
$\bullet$ For $d=3$, the constants in \eqref{estduran2} admit the following upper bounds (see Section \ref{apB}):
$$
\begin{aligned}
& A_{11} = A_{22} = A_{33} < 7.57, \qquad A_{12} = A_{13} = A_{21} = A_{23} = A_{31} = A_{32} < 3.5, \\[3pt]
& \widetilde{A}_{11} = \widetilde{A}_{22} = \widetilde{A}_{33} < \dfrac{1.21}{r^{3/2}}, \qquad \widetilde{A}_{12} = \widetilde{A}_{13} = \widetilde{A}_{21} = \widetilde{A}_{23} = \widetilde{A}_{31} = \widetilde{A}_{32} < \dfrac{0.68}{r^{3/2}},\\[3pt]
& B_{11} = B_{21} = B_{31} = B_{12} = B_{22} = B_{32} = B_{13} = B_{23} = B_{33} < \dfrac{11.2}{r},\\[3pt]
& \widetilde{B}_{11} = \widetilde{B}_{21} = \widetilde{B}_{31} = \widetilde{B}_{12} = \widetilde{B}_{22} = \widetilde{B}_{32} = \widetilde{B}_{13} = \widetilde{B}_{23} = \widetilde{B}_{33} < \dfrac{1.97}{r^{5/2}}.
\end{aligned}
$$
By inserting these values into \eqref{estduran1} we obtain:
\begin{equation} \label{estduran33d}
\left\lbrace
\begin{aligned}
& \| T_{11,1}(g) \|_{L^{2}(\mathcal{O})} = \| T_{22,1}(g) \|_{L^{2}(\mathcal{O})} = \| T_{33,1}(g) \|_{L^{2}(\mathcal{O})} \leq \left( 15.14 + \dfrac{3.43}{r^{3/2}} \sqrt{| \mathcal{O} |} \right) \| g \|_{L^{2}(\mathcal{O})}, \\[6pt]
& \| T_{kj,1}(g) \|_{L^{2}(\mathcal{O})} \leq \left( 7 + \dfrac{1.93}{r^{3/2}} \sqrt{| \mathcal{O} |} \right) \| g \|_{L^{2}(\mathcal{O})} \qquad \text{for every} \ k,j \in \{1,2,3\}, \ k \neq j, \\[6pt]
& \| T_{kj,2}(g) \|_{L^{2}(\mathcal{O})} \leq \dfrac{\delta(\mathcal{O})}{r} \left( 22.4 + \dfrac{5.58}{r^{3/2}} \sqrt{| \mathcal{O} |} \right) \| g \|_{L^{2}(\mathcal{O})} \ \ \ \ \ \ \text{for every} \ k,j \in \{1,2,3\}.
\end{aligned}
\right.
\end{equation}
Finally, the conclusion is  obtained 
by  inserting the estimates  \eqref{estduran3}-\eqref{estduran33d} into \eqref{f:Young}.  
\qed
\\

\medskip
We are now in a position to give the

\vfill\eject
\noindent {\bf Proof of Theorem \ref{unionstarshapedcor}}.
We have to show that the claim in \eqref{f:bp} is fulfilled if one takes $C$ equal to the constant $M$ defined in the statement.
Thus, for a given $g \in L ^ 2 _ 0 (\Omega_{0})$, we are going to construct a vector field $v \in H ^ 1_0 (\Omega_{0})$ such that
$\nabla \cdot v = g$ in  $\Omega_{0} $, and $\| \nabla v \| _{ L ^ 2 (\Omega_{0}) } \leq M \| g \| _{L ^ 2 (\Omega_{0})}$. For the sake of clearness, we
divide the procedure into four steps.

\medskip
\noindent
{\it 1)  Domain decomposition}.
We write $\Omega_{0} = \Omega _1 \cup \Omega_2$, where $\Omega _1$ and $\Omega _2$ are star-shaped with respect to some ball.
Essentially, $\Omega_{1}$ is the region lying above $P$ which is ``illuminated'' by a ball placed in an upper corner of  $Q$, 
while $\Omega_{2}$ is the region lying below $P$ which is ``illuminated'' by a ball placed in the opposite corner. 
To give a precise analytic description of these sets, it is mandatory to distinguish the cases $d = 2$ and $d = 3$, 
having in mind Figure \ref{stars}.

\begin{itemize}[leftmargin=*]
\item{}  For $d=2$, we denote by $T(P_{1}; P_{2}; P_{3}) \subset \mathbb{R}^{2}$ the triangle with vertices $P_{1}$, $P_{2}$ and $P_{3}$.  We set
\begin{itemize} [leftmargin=*]
	\item[$\triangleright$] $\Omega_ {1} \doteq [(-L,-a) \times (-L,L)] \cup [(-a,L) \times (b,L)] \cup T ((a,b); (L,b); (L, b - \alpha_{*}) )$. This domain  is star-shaped with respect to the disk $\left( x + \frac{L+a}{2} \right)^2 + \left( y - \frac{L+a}{2} \right)^2 < \left( \frac{L-a}{2} \right)^2$, provided the point $(L, b - \alpha_{*})$ is defined as the intersection between the tangent line to such a disk through the vertex
	$(a, b)$  and the line $\{ x = L\}$.   Some lengthy computations show that
	$$ \alpha_{*} = \frac{L-a}{8a(L+a)} \left[ (L+3a)(L+a-2b) - (L-a) \sqrt{(L+a-2b)^{2} + 8a(L+a)} \right] \,.$$
Note in particular that $\alpha _* \geq 0$ since $a \geq b$.
		\item[$\triangleright$] $\Omega_ {2} \doteq [(a,L) \times (-L,L)] \cup [(-L,a) \times (-L,-b)] \cup T ((-a,-b); (-L,-b); (-L, \alpha_{*}-b ))$, with $\alpha_*$ defined as above. The domain $\Omega_ {2}$ is star-shaped with respect to the disk $\left( x - \frac{L+a}{2} \right)^2 + \left( y + \frac{L+a}{2} \right)^2 < \left( \frac{L-a}{2} \right)^2$.
\end{itemize}
\noindent
Some tedious computations give
\begin{equation} \label{sigmarea}
| \Omega_{1} | = | \Omega_{2} |=\sigma_{2} \, , \qquad | \Omega_{1} \cap \Omega_{2} | =  \gamma_{2} \,,
\end{equation}
with $\sigma_{2}$ and $\gamma_{2}$ defined as in \eqref{sigmagamma}.

\item{}
For $d= 3$, in order to avoid too lengthy computations, we opt for a simpler decomposition: we set
\begin{itemize} [leftmargin=*]
	\item[$\triangleright$] $\Omega_ {1} \doteq [(-L,L) \times (-L,L) \times (c,L)] \cup [(a,L) \times (-L,L) \times (-L,c)] \cup [(-L,a) \times (b,L) \times (-L,c)]$. This domain is star-shaped with respect to the ball
	$$
	\left( x - \dfrac{L+a}{2} \right)^2 + \left( y - \dfrac{L+b}{2} \right)^2 + \left( z - \dfrac{L+c}{2} \right)^2 < \left( \dfrac{L-a}{2} \right)^2.
	$$
	\item[$\triangleright$] $\Omega_ {2} \doteq [(-L,L) \times (-L,L) \times (-L,-c)] \cup [(-L,-a) \times (-L,L) \times (-c,L)] \cup [(-a,L) \times (-L,-b) \times (-c,L)]$. The domain $\Omega _2$ is star-shaped with respect to the ball
	$$
	\left( x + \dfrac{L+a}{2} \right)^2 + \left( y + \dfrac{L+b}{2} \right)^2 + \left( z + \dfrac{L+c}{2} \right)^2 < \left( \dfrac{L-a}{2} \right)^2.
	$$
\end{itemize}
\noindent
In this case, again via direct computations, we obtain 
\begin{equation} \label{sigmarea3d}
| \Omega_{1} | = | \Omega_{2} |=\sigma_{3} \, , \qquad | \Omega_{1} \cap \Omega_{2} | =  \gamma_{3} \,.
\end{equation}
with $\sigma_3$ and $\gamma_3$ as in \eqref{sigmagamma3d}.
\end{itemize}

\noindent
{\it 2) Decomposition of the datum $g$}.
We argue as in \cite{bogovskii1980}, see also \cite[Lemma III.3.2 and Theorem III.3.1]{galdi2011introduction},
 and we decompose $g$ as
$$
g = g_1 + g_2 \ \ \text{in} \ \ \Omega_{0}, \ \ \ \ g_{1} \in L_{0}^{2}(\Omega_{1}) \ \text{with} \ \text{supp}(g_{1}) \subset \overline{\Omega_{1}}, \ \ \ \ g_{2} \in L_{0}^{2}(\Omega_{2}) \ \text{with} \ \text{supp}(g_{2}) \subset \overline{\Omega_{2}}.
$$
The functions $g_{1}, g_ {2} : \Omega_{0} \longrightarrow \mathbb{R}$ are explicitly defined by
$$
g_{1}(\xi) =
\left\lbrace
\begin{array}{ll}
g(\xi) - \dfrac{\chi^{*}(\xi)}{| \Omega_{1} \cap \Omega_{2} |} \displaystyle\int_{\Omega_{1}} g(\xi') \,d\xi' \ \ & \text{if} \ \xi \in \Omega_{1}, \\
0 \ \ & \text{if} \ \xi \in \Omega_{2} \setminus \Omega_{1},
\end{array}
\right.
$$
\vspace{2mm}
$$
g_{2}(\xi) =
\left\lbrace
\begin{array}{ll}
[1 - \chi^{*}(\xi)]g(\xi) - \dfrac{\chi^{*}(\xi)}{| \Omega_{1} \cap \Omega_{2} |} \displaystyle\int_{\Omega_{2} \setminus \Omega_{1}} g(\xi') \,d\xi' \ \ & \text{if} \ \xi \in \Omega_{2}, \\
0 \ \ & \text{if} \ \xi \in \Omega_{1} \setminus \Omega_{2},
\end{array}
\right.
$$
where $\chi^{*}$ is the characteristic function of the set $\Omega_{1} \cap \Omega_{2}$. We have
\begin{equation} \label{decompgaldi4}
\| g_{1} \|_{L^{2}(\Omega_{0})} = \alpha_g\, ,\qquad\| g_{2} \|_{L^{2}(\Omega_{0})} = \beta_g\, ,
\end{equation}
where
\begin{equation} \label{normas}
\begin{aligned}
& \alpha_{g} \doteq \left\{ \| g \|^{2}_{L^{2}(\Omega_{1})} + \dfrac{1}{| \Omega_{1} \cap \Omega_{2} |} \left( \int_{\Omega_{1}} g(\xi)\,d\xi \right) \left[\int_{\Omega_{1}} g(\xi)\,d\xi - 2 \int_{\Omega_{1} \cap \Omega_{2}} g(\xi)\,d\xi \right] \right\}^{1/2} , \\[7pt]
& \beta_{g} \doteq \left\{ \| g \|^{2}_{L^{2}(\Omega_{2} \setminus \Omega_{1})} + \dfrac{1}{| \Omega_{1} \cap \Omega_{2} |} \left( \int_{\Omega_{2} \setminus \Omega_{1}} g(\xi)\,d\xi \right)^{2} \right\}^{1/2}.
\end{aligned}
\end{equation}
Notice that for every $g \in L_{0}^{2}(\Omega_{0})$, we have
\begin{equation} \label{normasprop}
\int_{\Omega_{1}} g(\xi)\,d\xi + \int_{\Omega_{2} \setminus \Omega_{1}} g(\xi)\,d\xi = 0.
\end{equation}	
In view of \eqref{normasprop},  and applying Jensen inequality, we get
$$
\begin{aligned}
\alpha_{g}^{2} + \beta_{g}^{2} & = \| g \|^{2}_{L^{2}(\Omega_{0})} + \dfrac{2}{\gamma_{d}} \left( \int_{\Omega_{1}} g(\xi)\,d\xi \right) \left( \int_{\Omega_{1} \setminus \Omega_{2}} g(\xi)\,d\xi \right) \leq \| g \|^{2}_{L^{2}(\Omega_{0})} + \dfrac{2}{\gamma} \left( \int_{\Omega_{0}} |g(\xi)| \,d\xi \right)^2 \\[6pt]
& \leq \left( 1 + \dfrac{2}{\gamma_{d}} |\Omega_{0}| \right) \, \| g \|^{2}_{L^{2}(\Omega_{0})} \, ,
\end{aligned}
$$
so that
\begin{equation} \label{jensen1}
\alpha_{g} + \beta_{g} \leq \sqrt{  2  \left( \alpha _ g ^ 2 + \beta_g ^ 2 \right)} \leq \sqrt{2 \left( 1 + \dfrac{2}{\gamma_{d}} |\Omega_{0}| \right)} \, \| g \|_{L^{2}(\Omega_{0})}.
\end{equation}
\medskip

\noindent
{\it 3) Solving two distinct Bogovskii problems}.
We deal with the Bogovskii problem  on each of the two domains $\Omega_1$ and $\Omega_2$.
\begin{itemize} [leftmargin=*]
\item{} For $d=2$, we have
$r=\frac{L-a}2$ and
$
\delta(\Omega_{1})=\delta(\Omega_{2})=2\sqrt{2}\, L$.
By  Proposition \ref{p:starshaped},  there exist two vector fields $v_{1} \in H_{0}^{1}(\Omega_{1})$ and $v_{2} \in H_{0}^{1}(\Omega_{2})$ verifying, for $k = 1, 2$,
$$
\nabla \cdot v_{k} = g_{k} \quad \text{in} \quad \Omega_{k} \, ,$$
$$
\| \nabla v_{k} \|_{L^{2}(\Omega_{k})} \leq 2 \Bigg[ 129.35 + \dfrac{143.86 \sqrt{\sigma_{2}}}{L-a} + \dfrac{45.36\sigma_{2}}{(L-a)^2}
 + \dfrac{64L^{2}}{(L-a)^2} \left( 13.79 + \dfrac{7.28 \sqrt{\sigma_{2}}}{L-a}  \right)^{2} \Bigg]^{1/2} \| g_{k} \|_{L^{2}(\Omega_{k})}\,.
$$
\item{} For $d= 3$,  we have $r=\frac{L-a}2$ and
$ \delta(\Omega_{1})=\delta(\Omega_{2})= 2\sqrt{3}\, L \,.$  By  Proposition \ref{p:starshaped},  there exist two vector fields $v_{1} \in H_{0}^{1}(\Omega_{1})$ and $v_{2} \in H_{0}^{1}(\Omega_{2})$ verifying, for $k = 1, 2$,
$$
\nabla \cdot v_{k} = g_{k} \quad \text{in} \quad \Omega_{k} \, ,$$
{\small
$$
\| \nabla v_{k} \|_{L^{2}(\Omega_{k})} \leq \sqrt{6} \Bigg[ 327.23 + \dfrac{445.17 \sqrt{\sigma_{3}}}{(L-a)^{3/2}} + \dfrac{153.85 \sigma_{3}}{(L-a)^3}
+ \dfrac{144 L^{2}}{(L-a)^2} \left( 22.4 + \dfrac{15.79 \sqrt{\sigma_{3}}}{(L-a)^{3/2}}  \right)^{2} \Bigg]^{1/2} \| g_{k} \|_{L^{2}(\Omega_{k})}\,.
$$}
\end{itemize}
\noindent
{\it 4) Glueing and conclusion}.  After extending both the fields $v_{1}$ and $v_{2}$ defined in Step 3 to zero, respectively outside $\Omega_{1}$ and $\Omega_{2}$,  we infer that the vector field $v_0\doteq v_1+v_2\in H_{0}^{1}(\Omega_{0})$, satisfies $\nabla \cdot v _0 = g$ in $\Omega_{0}$, along with  the following bounds
\begin{itemize}[leftmargin=*]
\item{} For $d= 2$,
\begin{equation} \label{final1}
\| \nabla v_0 \|_{L^{2}(\Omega_{0})} \leq 2 \Bigg[ 129.35 + \dfrac{143.86 \sqrt{\sigma_{2}}}{L-a} + \dfrac{45.36\sigma_{2}}{(L-a)^2}
+ \dfrac{64L^{2}}{(L-a)^2} \left( 13.79 + \dfrac{7.28 \sqrt{\sigma_{2}}}{L-a} \right)^{2} \Bigg]^{1/2} (\alpha_{g} + \beta_{g}).
\end{equation}
\item{} For $d = 3$,
{\small
\begin{equation} \label{final3d}
\| \nabla v_{0} \|_{L^{2}(\Omega_{0})} \leq \sqrt{6} \Bigg[ 327.23 + \dfrac{445.17 \sqrt{\sigma_{3}}}{(L-a)^{3/2}} + \dfrac{153.85 \sigma_{3}}{(L-a)^3}
+ \dfrac{144 L^{2}}{(L-a)^2} \left( 22.4 + \dfrac{15.79 \sqrt{\sigma_{3}}}{(L-a)^{3/2}}  \right)^{2} \Bigg]^{1/2} (\alpha_{g} + \beta_{g}).
\end{equation}}
\end{itemize}
The conclusion then follows by inserting \eqref{jensen1} into \eqref{final1} and \eqref{final3d}.

\section{Estimates for the norms of mollifiers} \label{normsmollifiers}

In this section we provide the  computations leading to the estimates of the norms appearing in \eqref{estduran2}, that we have used in the proof of Proposition \ref{p:starshaped}. 
We observe that the radial mollifier introduced in \eqref{mollifier}-\eqref{mollifier3d} can be rewritten as 
$$
\omega_{d}(\xi) = \dfrac{\ell_{d}}{r^{d}} \, \omega_{0}\left( \dfrac{\xi}{r} \right) \qquad \forall \xi \in \mathbb{R}^d,
$$
where $\ell_{d} > 0$ is the normalization constant in \eqref{normalcons}, 
and $\omega_{0} \in \mathcal{C}^{\infty}_{0}(\mathbb{R}^d)$ is the standard radial mollifier
supported in the unit ball $B_{0}$ of $\mathbb{R}^d$, namely
$$
\omega_{0}(\xi) =
\left\lbrace
\begin{aligned}
& \exp\left( \dfrac{1}{| \xi |^{2} - 1} \right) \ \ \ \ & \text{if} \ \ | \xi | < 1, \\[6pt]
& 0 & \text{if} \ \ | \xi | \geq 1.
\end{aligned}
\right.
$$
Thus, in order to compute all the norms appearing in \eqref{estduran2}, after writing the corresponding integrals over $B(0,r) \subset \mathbb{R}^d$, via   
the change of variables $ \zeta = \xi / r $  we reduce ourselves to  integrals over the unit ball $B_{0} \subset \mathbb{R}^d$. These integrals are independent of $r$, and they can be easily computed numerically, yielding the following bounds. 

\subsection{2d case}\label{apA}
Let $\omega_{2}$ be the mollifier defined in \eq{mollifier}, with $\ell_{2} \approx 2.14357$.\par\noindent
$\bullet$ \textit{Bounds for the norms of zeroth-order derivatives:}
$$
\begin{aligned}
& \| \omega_{2} \|_{L^{1}(\mathcal{B})} = 1; \ \ \ \ \ \| x\, \omega_{2} \|_{L^{1}(\mathcal{B})} = \| y\, \omega_{2} \|_{L^{1}(\mathcal{B})} = r \ell_{2} \| x\, \omega_{0} \|_{L^{1}(B_{0})} < 0.31r.
\end{aligned}
$$
\noindent
$\bullet$ \textit{Bounds for the norms of first-order derivatives:}
$$
\begin{aligned}
& \left\| \dfrac{\partial \omega_{2}}{\partial x} \right\|_{L^{1}(\mathcal{B})} = \left\| \dfrac{\partial \omega_{2}}{\partial y} \right\|_{L^{1}(\mathcal{B})} = \dfrac{\ell_{2}}{r} \left\| \dfrac{\partial \omega_{0}}{\partial x} \right\|_{L^{1}(B_{0})} < \dfrac{1.91}{r};\\[6pt]
& \left\| \dfrac{\partial \omega_{2}}{\partial x} \right\|_{L^{\infty}(\mathcal{B})} = \left\| \dfrac{\partial \omega_{2}}{\partial y} \right\|_{L^{\infty}(\mathcal{B})} = \dfrac{\ell_{2}}{r^3} \left\| \dfrac{\partial \omega_{0}}{\partial x} \right\|_{L^{\infty}(B_0)}  < \dfrac{1.72}{r^3}; \\[6pt]
& \left\| \dfrac{\partial}{\partial x} (x \omega_{2}) \right\|_{L^{1}(\mathcal{B})} = \left\| \dfrac{\partial}{\partial y} (y \omega_{2}) \right\|_{L^{1}(\mathcal{B})} = \ell_{2} \left\| \dfrac{\partial}{\partial x} (x \omega_{0}) \right\|_{L^{1}(B_{0})} < 1.18; \\[6pt]
\end{aligned}
$$
$$
\begin{aligned}
& \left\| \dfrac{\partial}{\partial x} (x \omega_{2}) \right\|_{L^{\infty}(\mathcal{B})} = \left\| \dfrac{\partial}{\partial y} (y \omega_{2}) \right\|_{L^{\infty}(\mathcal{B})} = \dfrac{\ell_{2}}{r^2} \left\| \dfrac{\partial}{\partial x} (x \omega_{0}) \right\|_{L^{\infty}(B_{0})} < \dfrac{1.19}{r^2};\\[6pt]
& \left\| \dfrac{\partial}{\partial x} (y \omega_{2}) \right\|_{L^{1}(\mathcal{B})} = \left\| \dfrac{\partial}{\partial y} (x \omega_{2}) \right\|_{L^{1}(\mathcal{B})} = \ell_{2} \left\| \dfrac{\partial}{\partial x} (y \omega_{0}) \right\|_{L^{1}(B_{0})} < 0.64;\\[6pt]
& \left\| \dfrac{\partial}{\partial x} (y \omega_{2}) \right\|_{L^{\infty}(\mathcal{B})} = \left\| \dfrac{\partial}{\partial y} (x \omega_{2}) \right\|_{L^{\infty}(\mathcal{B})} =  \dfrac{\ell_{2}}{r^{2}} \left\| \dfrac{\partial}{\partial x} (y \omega_{0}) \right\|_{L^{\infty}(B_{0})} < \dfrac{0.67}{r^2}.\\[6pt]
\end{aligned}
$$
\medskip\noindent
$\bullet$ \textit{Bounds for the norms of second-order derivatives:}
$$
\begin{aligned}
& \left\| \dfrac{\partial^2 \omega_{2}}{\partial x^2} \right\|_{L^{1}(\mathcal{B})} = \left\| \dfrac{\partial^2 \omega_{2}}{\partial y^2} \right\|_{L^{1}(\mathcal{B})} = \dfrac{\ell_{2}}{r^2} \left\| \dfrac{\partial^2 \omega_{0}}{\partial x^2} \right\|_{L^{1}(B_{0})} < \dfrac{8.75}{r^2};\\[6pt]
& \left\| \dfrac{\partial^{2} }{\partial x^{2}}(x \omega_{2}) \right\|_{L^{1}(\mathcal{B})} = \left\| \dfrac{\partial^{2} }{\partial y^{2}}(y \omega_{2}) \right\|_{L^{1}(\mathcal{B})} = \dfrac{\ell_{2}}{r} \left\| \dfrac{\partial^{2} }{\partial x^{2}}(x \omega_{0}) \right\|_{L^{1}(B_{0})} < \dfrac{6.98}{r};\\[6pt]
& \left\| \dfrac{\partial^{2} }{\partial x^{2}}(y \omega_{2}) \right\|_{L^{1}(\mathcal{B})} = \left\| \dfrac{\partial^{2} }{\partial y^{2}}(x \omega_{2}) \right\|_{L^{1}(\mathcal{B})} = \dfrac{\ell_{2}}{r} \left\| \dfrac{\partial^{2} }{\partial x^{2}}(y \omega_{0}) \right\|_{L^{1}(B_{0})} < \dfrac{3.08}{r}.
\end{aligned}
$$

\subsection{3d case}\label{apB}
Let $\omega_3$ be the mollifier defined in \eq{mollifier3d}, with $\ell_{3} \approx 2.26712$.\par\noindent
$\bullet$ \textit{Bounds for the norms of zeroth-order derivatives:}
$$
\begin{aligned}
& \| \omega_{3} \|_{L^{1}(\mathcal{B})} = 1; \ \ \ \ \ \| x\, \omega_{3} \|_{L^{1}(\mathcal{B})} = \| y\, \omega_{3} \|_{L^{1}(\mathcal{B})} = \| z\, \omega_{3} \|_{L^{1}(\mathcal{B})}  = r \ell_{3} \| x\, \omega_{0} \|_{L^{1}(B_{0})} < 0.28r.
\end{aligned}
$$
\noindent
$\bullet$ \textit{Bounds for the norms of first-order derivatives:}
$$
\begin{aligned}
& \left\| \dfrac{\partial \omega_{3}}{\partial x} \right\|_{L^{1}(\mathcal{B})} =  \left\| \dfrac{\partial \omega_{3}}{\partial y} \right\|_{L^{1}(\mathcal{B})} = \left\| \dfrac{\partial \omega_{3}}{\partial z} \right\|_{L^{1}(\mathcal{B})} = \dfrac{\ell_{3}}{r} \left\| \dfrac{\partial \omega_{0}}{\partial x} \right\|_{L^{1}(B_{0})} < \dfrac{2.12}{r};\\[6pt]
& \left\| \dfrac{\partial \omega_{3}}{\partial x} \right\|_{L^{\infty}(\mathcal{B})} = \left\| \dfrac{\partial \omega_{3}}{\partial y} \right\|_{L^{\infty}(\mathcal{B})} = \left\| \dfrac{\partial \omega_{3}}{\partial z} \right\|_{L^{\infty}(\mathcal{B})} = \dfrac{\ell_{3}}{r^4} \left\| \dfrac{\partial \omega_{0}}{\partial x} \right\|_{L^{\infty}(B_0)} < \dfrac{1.82}{r^4};
\end{aligned}
$$
$$
\begin{aligned}
& \left\| \dfrac{\partial}{\partial x} (x \omega_{3}) \right\|_{L^{1}(\mathcal{B})} = \left\| \dfrac{\partial}{\partial y} (y \omega_{3}) \right\|_{L^{1}(\mathcal{B})} = \left\| \dfrac{\partial}{\partial z} (z \omega_{3}) \right\|_{L^{1}(\mathcal{B})} = \ell_{3} \left\| \dfrac{\partial}{\partial x} (x \omega_{0}) \right\|_{L^{1}(B_{0})} < 1.16;\\[6pt]
& \left\| \dfrac{\partial}{\partial x} (x \omega_{3}) \right\|_{L^{\infty}(\mathcal{B})} = \left\| \dfrac{\partial}{\partial y} (y \omega_{3}) \right\|_{L^{\infty}(\mathcal{B})} = \left\| \dfrac{\partial}{\partial z} (z \omega_{3}) \right\|_{L^{\infty}(\mathcal{B})} = \dfrac{\ell_{3}}{r^3} \left\| \dfrac{\partial}{\partial x} (x \omega_{0}) \right\|_{L^{\infty}(B_{0})} < \dfrac{1.26}{r^3};
\end{aligned}
$$
$$
\begin{aligned}
& \left\| \dfrac{\partial}{\partial x} (y \omega_{3}) \right\|_{L^{1}(\mathcal{B})} = \left\| \dfrac{\partial}{\partial x} (z \omega_{3}) \right\|_{L^{1}(\mathcal{B})} = \left\| \dfrac{\partial}{\partial y} (x \omega_{3}) \right\|_{L^{1}(\mathcal{B})} = \ell_{3} \left\| \dfrac{\partial}{\partial x} (y \omega_{0}) \right\|_{L^{1}(B_{0})} < 0.64;\\[6pt]
& \left\| \dfrac{\partial}{\partial y} (z \omega_{3}) \right\|_{L^{1}(\mathcal{B})} = \left\| \dfrac{\partial}{\partial z} (x \omega_{3}) \right\|_{L^{1}(\mathcal{B})} = \left\| \dfrac{\partial}{\partial z} (y \omega_{3}) \right\|_{L^{1}(\mathcal{B})} = \ell_{3} \left\| \dfrac{\partial}{\partial x} (y \omega_{0}) \right\|_{L^{1}(B_{0})} < 0.64;\\[6pt]
& \left\| \dfrac{\partial}{\partial x} (y \omega_{3}) \right\|_{L^{\infty}(\mathcal{B})} = \left\| \dfrac{\partial}{\partial x} (z \omega_{3}) \right\|_{L^{\infty}(\mathcal{B})} = \left\| \dfrac{\partial}{\partial y} (x \omega_{3}) \right\|_{L^{\infty}(\mathcal{B})} = \dfrac{\ell_{3}}{r^{3}} \left\| \dfrac{\partial}{\partial x} (y \omega_{0}) \right\|_{L^{\infty}(B_{0})} < \dfrac{0.71}{r^{3}};\\[6pt]
& \left\| \dfrac{\partial}{\partial y} (z \omega_{3}) \right\|_{L^{\infty}(\mathcal{B})} = \left\| \dfrac{\partial}{\partial z} (x \omega_{3}) \right\|_{L^{\infty}(\mathcal{B})} = \left\| \dfrac{\partial}{\partial z} (y \omega_{3}) \right\|_{L^{\infty}(\mathcal{B})} = \dfrac{\ell_{3}}{r^{3}} \left\| \dfrac{\partial}{\partial x} (y \omega_{0}) \right\|_{L^{\infty}(B_{0})} < \dfrac{0.71}{r^{3}};\\[6pt]
\end{aligned}
$$
\medskip\noindent
$\bullet$ \textit{Bounds for the norms of second-order derivatives:}
$$
\begin{aligned}
& \left\| \dfrac{\partial^2 \omega_{3}}{\partial x^2} \right\|_{L^{1}(\mathcal{B})} = \left\| \dfrac{\partial^2 \omega_{3}}{\partial y^2} \right\|_{L^{1}(\mathcal{B})} = \left\| \dfrac{\partial^2 \omega_{3}}{\partial z^2} \right\|_{L^{1}(\mathcal{B})} = \dfrac{\ell_{3}}{r^2} \left\| \dfrac{\partial^2 \omega_{0}}{\partial x^2} \right\|_{L^{1}(B_{0})} < \dfrac{10.22}{r^2};\\[6pt]
& \left\| \dfrac{\partial^{2} }{\partial x^{2}}(x \omega_{3}) \right\|_{L^{1}(\mathcal{B})} = \left\| \dfrac{\partial^{2} }{\partial y^{2}}(y \omega_{3}) \right\|_{L^{1}(\mathcal{B})} = \left\| \dfrac{\partial^{2} }{\partial z^{2}}(z \omega_{3}) \right\|_{L^{1}(\mathcal{B})} = \dfrac{\ell_{3}}{r} \left\| \dfrac{\partial^{2} }{\partial x^{2}}(x \omega_{0}) \right\|_{L^{1}(B_{0})} < \dfrac{7.29}{r};\\[6pt]
\end{aligned}
$$
$$
\begin{aligned}
& \left\| \dfrac{\partial^{2} }{\partial x^{2}}(y \omega_{3}) \right\|_{L^{1}(\mathcal{B})} = \left\| \dfrac{\partial^{2} }{\partial x^{2}}(z \omega_{3}) \right\|_{L^{1}(\mathcal{B})} = \left\| \dfrac{\partial^{2} }{\partial y^{2}}(x \omega_{3}) \right\|_{L^{1}(\mathcal{B})} = \dfrac{\ell_{3}}{r} \left\| \dfrac{\partial^{2} }{\partial x^{2}}(y \omega_{0}) \right\|_{L^{1}(B_{0})} < \dfrac{3.22}{r}; \\[6pt]
& \left\| \dfrac{\partial^{2} }{\partial y^{2}}(z \omega_{3}) \right\|_{L^{1}(\mathcal{B})} = \left\| \dfrac{\partial^{2} }{\partial z^{2}}(x \omega_{3}) \right\|_{L^{1}(\mathcal{B})} = \left\| \dfrac{\partial^{2} }{\partial z^{2}}(y \omega_{3}) \right\|_{L^{1}(\mathcal{B})} = \dfrac{\ell_{3}}{r} \left\| \dfrac{\partial^{2} }{\partial x^{2}}(y \omega_{0}) \right\|_{L^{1}(B_{0})} < \dfrac{3.22}{r}.
\end{aligned}
$$

\section{Proof of Theorem \ref{teononconstant}} \label{finalproof}

We proceed as outlined in Section \ref{outline}. Firstly we denote by $A_1  \in H^{1}(Q) \cap \mathcal{C}(\overline{Q})$ the extension of $h$ to $Q$ given by Theorem \ref{gadliardo}, so that the norm of $A_1$ in $H ^ 1 (Q)$ can be explicitly computed in terms of the restrictions of $h$ to the $(d-1)$-dimensional faces of $Q$. Secondly, we take  the function $\phi \in W^{1,\infty}(Q)$ defined as in Section \ref{sec:phi} (by \eqref{def:phi2}-\eqref{def:phi3}), 
which satisfies in particular
$\phi = 1$ on $\partial Q$ and $\phi = 0$  in $\overline P$.
Recall from  \eqref{f:phi2} and \eqref{f:phi3} that, both for $d=2$ and $d= 3$, 
\begin{equation}\label{f:norme}
\| \phi \| _{L ^ \infty (Q)}  = 1 \qquad \text{ and } \qquad \| \nabla \phi \| _{L ^ \infty (Q)}  =
\frac{1}{L-a}
\,.
\end{equation}
The vector field  $A_{2} \doteq \phi A_{1}$,  satisfies
$A_{2}=h$ on $\partial Q$ and $A_2 = 0$  in $\overline P$. Moreover, $\nabla \cdot A_ 2  \in L^ 2_0 (\Omega _0)$, since \eqref{zeroflux2}   implies
$$\int _{\Omega_0} \nabla \cdot A_ 2  = \int _{\partial \Omega _0} A_2 \cdot \hat n = \int _{\partial Q}  h \cdot \hat n = 0\,.$$
Thirdly, let $A_3 \in H_{0}^{1}(\Omega_{0})$ be the solution to Bogovskii's problem with right-hand side given by $-\nabla \cdot A_2$:
\begin{equation} \label{nonconstant1}
\nabla \cdot A_{3} = - \nabla \cdot A_2 \qquad \text{and} \qquad \| \nabla A_{3} \|_{L^{2}(\Omega_{0})} \leq M \, \| \nabla \cdot A _2 \|_{L^{2}(\Omega_{0})},
\end{equation}
where $M$ is given by Proposition \ref{unionstarshapedcor}. Since
$$
\nabla \cdot A _ 2 = \nabla \cdot (\phi A _1) = \phi (\nabla \cdot A_{1}) + \nabla \phi \cdot A_{1} \quad \text{in} \quad \Omega_{0},
$$
we have
\begin{equation} \label{nonconstant2}
\begin{aligned}
\| \nabla \cdot A_2 \|_{L^{2}(\Omega_{0})} & \leq \| \phi  (\nabla \cdot A_{1}) \|_{L^{2}(\Omega_{0})} + \| \nabla \phi \cdot A_{1} \|_{L^{2}(\Omega_{0})} \\[6pt]
& \leq \|\phi \|_{L^{\infty}(\Omega_{0})} \, \| \nabla \cdot A_{1} \|_{L^{2}(Q)} + \| \nabla \phi \|_{L^{\infty}(\Omega_{0})} \, \| A_{1} \|_{L^{2}(Q)}
\\[6pt]
\mbox{(by \eqref{f:norme}) } & =\| \nabla \cdot A_{1} \|_{L^{2}(Q)} + \frac{1}{L-a} \, \| A_{1} \|_{L^{2}(Q)}\,.
\end{aligned}
\end{equation}
We point out that the function  $ \phi$  in \eqref{def:phi2} or \eqref{def:phi3}  minimizes the second line in \eqref{nonconstant2} and, in this respect, the choice of $\phi$ is optimal thanks to Theorem \ref{thm:infty}.
By construction, the function $v_0 \doteq A_2 + A_3$ satisfies
$$
v_0 = h  \ \text{ on } \partial Q\,, \quad v= 0 \ \text{ on } \partial P \,, \quad \nabla \cdot v _0 = 0 \ \text{ in } \Omega _0\, , \quad \nabla v_{0} = \nabla \phi \otimes A_1 + \phi  \nabla A_1 + \nabla A_3 \ \text{ in } \Omega_{0} \, .
$$
From \eqref{nonconstant1} and \eqref{nonconstant2} we derive the estimate
\begin{equation} \label{nonconstant3}
\begin{aligned}
\| \nabla v_{0} \|_{L^{2}(\Omega_{0})} & \leq \| \nabla \phi \otimes A_1 \|_{L^{2}(\Omega_{0})} + \| \phi \nabla A_1 \|_{L^{2}(\Omega_{0})} + \| \nabla A_3 \|_{L^{2}(\Omega_{0})} \\[6pt]
& \leq \frac{1}{L-a}  \, \| A_{1} \|_{L^{2}(Q)} +  \| \nabla A_{1} \|_{L^{2}(Q)} + M  \Big(  \| \nabla \cdot A_{1} \|_{L^{2}(Q)} + \frac{1}{L-a} \, \| A_{1} \|_{L^{2}(Q)} \Big).
\end{aligned}
\end{equation}
After extending it to zero to  $P \setminus \overline{K}$, $v_0$ becomes a solution to problem \eqref{div}.  
Indeed, it is clear that it matches the boundary condition $v_0= h$ on $\partial \Omega$, and that it is divergence free separately in $\Omega \setminus P$ and in $P \setminus K$. It is also readily checked that the equation $\nabla \cdot  v_0 = 0$ is satisfied in distributional sense in $\Omega$ since, for every test function $\psi \in \mathcal C ^ \infty _0 (\Omega)$, we have
$$\int_\Omega v _0 \cdot \nabla \psi = \int_{Q \setminus P} v _0 \cdot  \nabla \psi + \int_{P \setminus K}  v _0 \cdot \nabla \psi   =  
\int_{Q \setminus P} v _0 \cdot \nabla \psi  =-  \int_{Q \setminus P}  \psi  (\nabla \cdot  v _0 ) + \int _{\partial (Q \setminus P)} \psi ( v_0 \cdot \hat{n} ) = 0\,.$$

Finally, we observe that the inequality \eqref{nonconstant3} can be  re-written as:
\begin{equation} \label{nonconstant4}
\| \nabla v_{0} \|_{L^{2}(\Omega)} = \| \nabla v_{0} \|_{L^{2}(\Omega_{0})}  \leq  \frac{1 + M}{L-a} \,  \| A_{1} \|_{L^{2}( Q)}
+  \| \nabla A_{1} \|_{L^{2}(Q)} + M  \, \| \nabla \cdot A_{1} \|_{L^{2}(Q)} .
\end{equation}
\noindent
By inserting \eqref{f:phi2}-\eqref{f:phi3}  into \eqref{nonconstant4}, Theorem \ref{teononconstant} is proved. \qed

\bigskip

\par\medskip\noindent
{\bf Acknowledgements.} The second Author is supported by the PRIN project {\em Direct and inverse problems for partial
	differential equations: theoretical aspects and applications} and by the GNAMPA group of the INdAM. The third Author is
supported by the Primus Research Programme PRIMUS/19/SCI/01, by the program GJ17-01694Y of the Czech National Grant Agency GA\v{C}R, and by the
University Centre UNCE/SCI/023 of the Charles University in Prague.

\phantomsection
\addcontentsline{toc}{section}{References}
\bibliographystyle{abbrv}
\bibliography{references}
\vspace{5mm}

\begin{minipage}{100mm}
	Ilaria Fragalà and Filippo Gazzola\\
	Dipartimento di Matematica\\
	Politecnico di Milano\\
	Piazza Leonardo da Vinci 32\\
	20133 Milan - Italy\\
	{\bf E-mail}: ilaria.fragala@polimi.it - filippo.gazzola@polimi.it\\
    \\
	Gianmarco Sperone\\
	Department of Mathematical Analysis\\
	Charles University in Prague\\
	Sokolovská 83\\
	186 75 Prague - Czech Republic\\
	{\bf E-mail}: sperone@karlin.mff.cuni.cz
\end{minipage}
\end{document}